\documentclass[ejs]{imsart}
\RequirePackage[OT1]{fontenc}
\RequirePackage{amsthm,amsmath}
\RequirePackage[numbers]{natbib}
\RequirePackage[colorlinks,citecolor=blue,urlcolor=blue]{hyperref}
\RequirePackage{hypernat}
\usepackage{amssymb,graphicx,url,bbm,pictexwd}

\def\th{\theta}
\def\D{{\cal D}}
\def\IK{I\!\!K}

\def\a{\alpha}
\def\b{\beta}
\def\R{I\!\!R}
\def\dd{\Delta}

\def\d{\delta}

\def\D{{\cal D}}

\def\G{{\mathbb G}}
\def\H{{\mathbb H}}

\def\P{{\mathbb P}}

\def\labda1{\lambda_1}
\def\labda2{\lambda_2}
\def\e{\epsilon}
\def\f{\phi}

\def\s{\sigma}

\def\comment#1{\relax}

\def\=in{\mathop{\rm =}}

\def\eop{\hfill\mbox{$\Box$}\newline}

\arxiv{1109.6176}

\startlocaldefs
\numberwithin{equation}{section}
\theoremstyle{plain}
\newtheorem{theorem}{Theorem}[section]

\newtheorem{lemma}{Lemma}[section]
\newtheorem{remark}{Remark}[section]
\newtheorem{definition}{Definition}[section]

\endlocaldefs

\begin{document}

\begin{frontmatter}
\title{Estimators for the interval censoring problem}
\runtitle{Interval censoring}

\begin{aug}
\author{\fnms{Piet} \snm{Groeneboom}\corref{}\ead[label=e1]{P.Groeneboom@tudelft.nl}}
\and
\author{\fnms{Tom} \snm{Ketelaars}\ead[label=e2]{tom.ketelaars@credit-suisse.com}}

\address{Address of the first author\\
Delft Institute of Applied Mathematics\\
Mekelweg 4, 2628 CD Delft\\
The Netherlands\\
\printead{e1}\\
\url{http://dutiosc.twi.tudelft.nl/~pietg/}}

\address{Address of the second author\\
Credit Suisse\\
Attn.: Tom Ketelaars\\
11 Madison Avenue\\
New York, NY, 10010\\
\printead{e2}}

\runauthor{P.\ Groeneboom and T.\ Ketelaars}
\affiliation{Delft University of Technology}
\end{aug}

\begin{abstract}
We study three estimators for the interval censoring case 2 problem, a histogram-type estimator, proposed in \cite{lucien:99}, the maximum likelihood estimator (MLE) and the smoothed MLE, using a smoothing kernel. Our focus is on the asymptotic distribution of the estimators at a fixed point. The estimators are compared in a simulation study.
\end{abstract}



\end{frontmatter}

\section{Introduction}
\label{sec:intro}
Let $X_1,\ldots,X_n$ be a sample of unobservable random variables
from an unknown distribution function $F_0$ on the interval
$[0,1]$. More generally, we could take an arbitrary closed interval $[a,b]$ as support for the underlying
distribution, but for the purposes of the development of the theory, we can just as well take $[0,1]$, as
is also done in \cite{lucien:99}.

Suppose that one can observe
$n$ pairs
$(T_i,U_i),$ independent of $X_i,$ with a joint density function $h$ on the upper triangle of the unit square,
for which the sum of the marginal densities is bounded away from zero. Moreover,
\begin{equation}
\label{def_indicators}
\Delta_{i1}=1_{\{X_i\le T_i\}},\qquad
\dd_{i2}=1_{\{T_i<X_i\le U_i\}},\qquad\dd_{i3}=1-\dd_{i,1}-\dd_{i,1},
\end{equation}
provide the only information one has on the position of the random variables
$X_i$ with respect to the observation times $T_i$ and $U_i$. In
this set-up we want to estimate the unknown distribution function $F_0$, generating the ``unobservables" $X_i$.
This setting is known as {\it interval censoring, case 2.}

The model of {\it current status data}, also known as {\it interval censoring, case 1}, has been thoroughly studied, and has a theory which is considerably simpler than the theory for the interval censoring, case 2, model. In the current status model one only has one observation time $T_i$, corresponding to the unobservable $X_i$, and the only information we have about $X_i$ is whether $X_i$ is to the left or to the right of $T_i$.

Although the present paper mainly focuses on the case 2 model, we start by discussing the current status model, in order to put this paper into a more general context and to explain why the case 2 model is so much harder to study. In the current status model, the only observations which are available to us are the pairs
$$
(T_i,\dd_i),\qquad\dd_i=1_{\{X_i\le T_i\}},
$$
so we do not observe $X_i$ itself, but only its ``current status" $\dd_i$. 
The nonparametric maximumum likelihood estimator, commonly denoted by NPMLE or just MLE, maximizes the (partial) log likelihood
$$
\sum_{i=1}^n \left\{\dd_i\log F(T_i)+(1-\dd_i)\log\left(1- F(T_i)\right)\right\},
$$
where the maximization is over all distribution functions $F$.

The MLE can be found in one step by computing the left-continuous slope of the greatest convex minorant of the cusum diagram of the points $(0,0)$ and the points
\begin{equation}
\label{cusum_combined}
\left(i,\sum_{j\le i}\dd_{(j)}\right),\,i=1,\dots,n,
\end{equation}
using a notation, introduced in \cite{GrWe:92}. Here $\dd_{(j)}$ denotes the indicator corresponding to the $j$th order statistic $T_{(j)}$. The theory for this estimator is further developed in \cite{GrWe:92}, where also the (non-normal) pointwise limit distribution is derived and it is shown that the rate of convergence is $n^{-1/3}$.

In contrast, there is no such one-step algorithm for computing the MLE in the case 2 situation, where one wants to maximize
$$
\sum_{i=1}^n \left\{\dd_{i1}\log F(T_i)+\dd_{i2}\log\left\{F(U_i)- F(T_i)\right\}+\dd_{i3}\log\left(1- F(U_i)\right)\right\}.
$$
over distribution functions $F$. One has to take recourse to iterative algorithms, for example the iterative convex minorant algorithm, introduced in \cite{GrWe:92} and further developed in \cite{geurt:98}.
Moreover, the MLE can possibly achieve a faster local rate of convergence than in the current status model, depending on properties of the bivariate distribution of the observation times $(T_i,U_i)$.

In the so-called {\it non-separated case}, the density of the pair of observation times $(T_i,U_i)$ is positive on the diagonal, meaning that we can have arbitrarily small observation intervals $[T_i,U_i]$. For this situation, \cite{lucien:99} proposes a simple piecewise constant
estimator for $F_0,$ with the purpose of showing that in this situation an estimator
can be constructed that achieves the $(n\log n)^{-1/3}$ convergence rate, which is optimal in a minimax sense, both using a global loss function , and using a local loss function for the estimation at a fixed point. In the {\it separated case}, the observation times $T_i$ and $U_i$ cannot become arbitrarily close: in this case there exists an $\e>0$ so that $U_i-T_i>\e$ for each $i$. In this case the convergence rate of Birg\'e's estimator is $n^{-1/3}$ again, which is also the minimax rate for the current status model. For both situations we derive the asymptotic behavior of Birg\'e's estimator, and compare this with the behavior of the MLE in a simulation study. The simulations show a better behavior of the MLE, probably caused by the local adaptivity of the MLE.

A common complaint about the MLEs is that under the conditions for which the local asymptotic distribution
result is derived, other estimators can be suggested, which in fact attain a faster rate of convergence. Such estimators are discussed for the current status model in, e.g., \cite{piet_geurt_birgit:10}, \cite{birgit_geurt_piet:11} and \cite{piet:11}. We introduce a similar estimator below for the case 2 model below, the smoothed maximum likelihood estimator (SMLE). The smoothed MLE is defined by
\begin{align}
\label{SMLE}
\tilde F_n^{ML}(t)=\int \IK\left((t-u)/b_n\right)\,d\hat F_n(u),
\end{align}
where
$$
\IK(u)=\int_{-\infty}^u K(w)\,dw=\left\{\begin{array}{lll}
0\,&,\,u< -1\\
\displaystyle{\int_{-1}^u K(w)\,dw}&,\,u\in[-1,1],\\
1\,&,\,u>1,
\end{array}
\right.
$$
letting $K$ be a smooth symmetric kernel, with support $[-1,1]$, like the triweight kernel
$$
K(u)=\tfrac{35}{32}\left(1-u^2\right)^31_{[-1,1]}(u),
$$
and taking the bandwidth $b_n\asymp n^{-1/5}$.
Note that
\begin{align*}
\tilde f_n^{ML}(t)\stackrel{\mbox{\small def}}=\frac{d}{dt}\tilde F_n^{ML}(t)=\frac1{b_n}\int K\left((t-u)/b_n\right)\,d\hat F_n(u)
\end{align*}
is an estimate of the density $f_0$ of the underlying distribution function $F_0$.

Analogously to what has been proved for the current status model, we expect the smoothed MLE  to converge at (at least) rate $n^{-2/5}$ under appropriate regularity conditions. 
It is an attractive alternative to the MLE and histogram-type estimator of \cite{lucien:99}. We give a heuristic discussion on this in section \ref{sec:MSLE}. 
Just as in \cite{GeGr:97} and \cite{GeGr:99}, the asymptotic variance depends on the solution of an integral equation. The asymptotic expressions for the variance, obtained by solving these equations numerically, give a rather good fit with the actually observed variances, as shown in section \ref{sec:MSLE}. The SMLE can probably also be used for a two-sample test for interval censored data, analogous to the two-sample test for current status data, introduced in \cite{piet:11}. The MSE of the smoothed MLE is much smaller than that of Birg\'e's estimator or the MLE for smooth underlying distribution functions, as is illustrated in the sections on the simulations.

\begin{figure}[!ht]
\begin{center}
\includegraphics[scale=0.5]{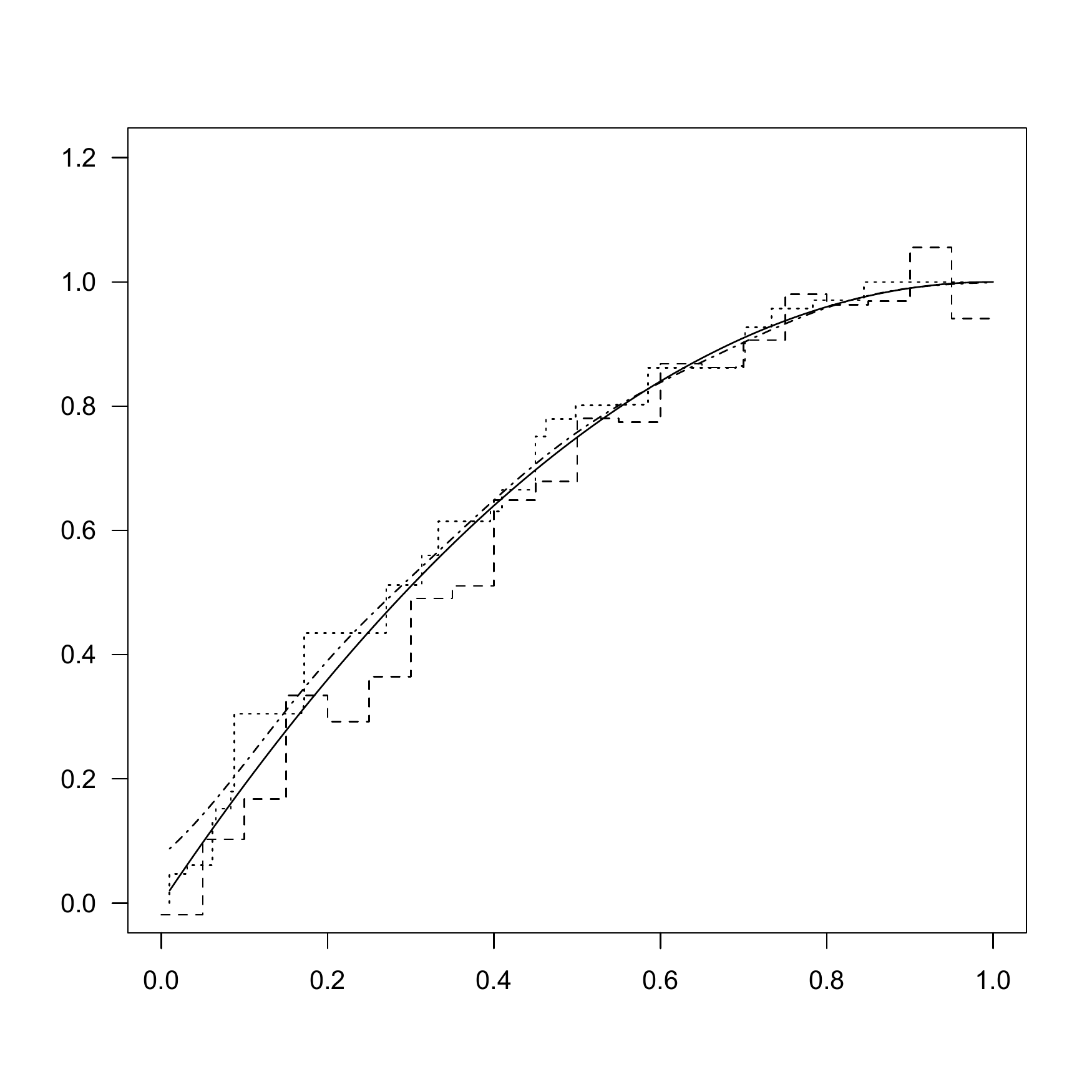}
\end{center}
\caption{Birg\'e's estimator (dashed), the MLE (dotted), and the smoothed MLE (dashed-dotted) for sample size $n=1000$ and $b_n=n^{-1/5}$, when $F_0(x)=1-(1-x)^2$ (solid curve) and the observation distribution is uniform on the upper triangle of the unit square.}
\label{fig:SMLE}
\end{figure}

A picture of the three estimators is shown in Figure \ref{fig:SMLE}. The MLE and smoothed MLE are monotone, in contrast with Birg\'e's
estimator. Also Birg\'e's estimator can have negative values and values larger than $1$; both events happen in the picture shown.
This cannot happen for the MLE and smoothed MLE, since these are based on isotonization; the smoothed MLE is an integral of
a positive kernel w.r.t.\ the (positive) jumps of the MLE, and inherits the monotonicity properties of the MLE. Although
histogram-type estimators (like Birg\'e's estimator) and kernel estimators without any isotonization are much easier to
analyze than the estimators, based on isotonization, the price one has to pay is the behavior illustrated in Figure \ref{fig:SMLE}.

\section{A local minimax result for the non-separated case}
\label{sec:minimax}
\setcounter{equation}{0}
In this section we derive a local minimax result for the
non-separated case of the interval censoring problem, case 2. This
result will provide the best possible local convergence rate and
also the best constant, as far as this constant depends on the underlying
distributions.

Our approach makes use of a perturbation $F_n$ of $F_0$ which is
defined by
$$
F_n(x)=\left\{ \begin{array}{ll}
F_0(x)&\textrm{if $x<t_0-c(n\log n)^{-1/3}$}\\
F_0(t_0-c(n\log n)^{-1/3}) &\textrm{if $x\in[t_0-c(n\log n)^{-1/3},t_0)$}\\
F_0(t_0+c(n\log n)^{-1/3}) &\textrm{if $x\in[t_0,t_0+c(n\log n)^{-1/3})$}\\
F_0(x)&\textrm{if $x\geq t_0+c(n\log n)^{-1/3}$}
\end{array} \right.
$$
for a $c>0$ to be specified below.

Before stating the theorem to be proved, we introduce some
notation. Let $\dd=(\dd_1,\dd_2)\in
\mathcal{T}:=\{(1,0),(0,1),(0,0)\}$ and define the densities $q_0$
and $q_n$ by
\begin{align*}
q_0(t,u,\underline{\delta})&=h(t,u)F_0(t)^{\delta_1}(F_0(u)-F_0(t))^{\delta_2}(1-F_0(u))^{1-\delta_1-\delta_2}\\
q_n(t,u,\underline{\delta})&=h(t,u)F_n(t)^{\delta_1}(F_n(u)-F_n(t))^{\delta_2}(1-F_n(u))^{1-\delta_1-\delta_2}
\end{align*}
with respect to the measure $\mu=\lambda_1 \otimes \lambda_2$ on
$\Omega=\mathbb{R}_{+}^{2} \times \mathcal{T},$ where $\lambda_1$
is the Lebesgue measure and $\lambda_2$ is counting measure. We
note that $q_0$ is the joint density of $(T,U,\dd_1,\dd_2)$.

Furthermore, let $(L_n),n\geq 1,$ be a sequence of estimators for
$F_0(t_0),$ based on samples of size n, generated by $q_0.$ That
is, we can write
$$
L_n=l_n((T_1,U_1,\Delta_{1,1},\Delta_{1,2}),\ldots,(T_n,U_n,\Delta_{n,1},\Delta_{n,1})),
$$
where $l_n$ is a Borel measurable function. Then, the following
theorem holds:
\begin{theorem}\label{mmax}
\begin{align*}
&\lim_{n\to \infty}\inf (n\log n)^{1/3}\max\{E_{n,q_0}|L_n-F_0(t_0)|,E_{n,q_n}|L_n-F_n(t_0)|\}\\
&\ge\frac{6^{1/3}}{4}\exp(-1/3)
\{f_0(t_0)^2/h(t_0,t_0)\}^{1/3},
\end{align*}
where $E_{n,q}$ denotes the expectation with respect to the
product measure $q^{\otimes n}.$
\end{theorem}
In our proof we need the following lemma, which is proved in
\cite{piet:96}. This type of result is often denoted as ``LeCam's lemma".

\begin{lemma}\label{mm}
Let $G$ be a set of probability densities on a measurable space
$(\Omega,A)$ with respect to a $\sigma$-finite dominating measure
$\mu$, and let L be a real-valued functional on $G.$ Moreover, let
$f:[0,\infty)\to \mathbb{R}$ be an increasing convex loss
function, with f(0)=0. Then, for any $q_1,q_2 \in G$ such that the
Hellinger distance $H(q_1,q_2)<1:$
\begin{align*}
&\inf_{L_n} \max \big\{E_{n,q_1} f(|L_n-Lq_1|),E_{n,q_2}f(|L_n-Lq_2|) \big\}\\
&\geq
f\left(\frac{1}{4}|Lq_1-Lq_2|\{1-H^2(q_1,q_2)\}^{2n}\right).
\end{align*}
\end{lemma}
\vspace*{0.5cm}

\begin{figure}[!ht]
 \label{fig:figure1}
\begin{center}
\strut\input{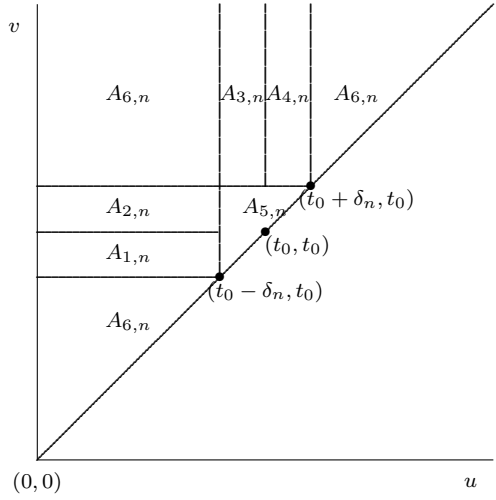} \caption{The areas
$A_{1,n},\ldots,A_{6,n}$}\end{center}
\end{figure}

\noindent
{\bf Proof} of theorem \ref{mmax}. Let the partitioning
$A_{1,n}\cup\ldots\cup A_{6,n}$ of $\{(t,u)\in
\mathbb{R}_{+}^{2}:t<u\}$ be defined by
\begin{align*}
A_{1,n}&=\{(t,u)\in \mathbb{R}_{+}^{2}:0<t<t_0-\delta_n,t_0-\delta_n\leq u <t_0)\}\\
A_{2,n}&=\{(t,u)\in \mathbb{R}_{+}^{2}:0<t<t_0-\delta_n,t_0\leq u < t_0+\delta_n\}\\
A_{3,n}&=\{(t,u)\in \mathbb{R}_{+}^{2}:t_0-\delta_n\leq t <t_0, t_0+\delta_n <u <\infty\}\\
A_{4,n}&=\{(t,u)\in \mathbb{R}_{+}^{2}:t_0\leq t <t_0+\delta_n, t_0+\delta_n <u <\infty\}\\
A_{5,n}&=\{(t,u)\in \mathbb{R}_{+}^{2}:t_0-\delta_n\leq t<t_0+\delta_n,t< u <t_0+\delta_n\}\\
A_{6,n}&=\{(t,u)\in\mathbb{R}_{+}^{2}:t<u\}\backslash\{A_{1,n}\cup
\ldots \cup A_{5,n}\},
\end{align*}
where $\delta_n=c(n \log n)^{-1/3}.$ The partitioning is shown in
figure  \ref{fig:figure1}.

Then the squared Hellinger distance between $q_0$ and $q_n$ can be
written as
\begin{align*}
H^2(q_n,q_0)&:=\frac{1}{2}\int_{\Omega}\left\{\sqrt{q_n}-\sqrt{q_0}\right\}d\mu\\
&= \frac{1}{2}\sum_{k=1}^{5}\int_{A_{k,n}}h(t,u) \left(\sqrt{F_n(t)}-\sqrt{F_0(t)}\right)^2 dt du+\\
&+\frac{1}{2}\sum_{k=1}^{5}\int_{A_{k,n}}h(t,u)\left(\sqrt{F_n(u)-F_n(t)}-\sqrt{F_0(u)-F_0(t)}\right)^2dt du\\
&+\frac{1}{2}\sum_{k=1}^{5}\int_{A_{k,n}}h(t,u)\left(\sqrt{1-F_n(u)}-\sqrt{1-F_0(u)}\right)^2dtdu.
\end{align*}
We now calculate the three integrals over $A_{1,n}.$

Obviously, we have
\begin{equation}\label{fin}
\int_{A_{1,n}}h(t,u)\left(\sqrt{F_n(t)}-\sqrt{F_0(t)}\right)^2 dt
du=0.
\end{equation}
Furthermore,
\begin{align*}
&\int_{A_{1,n}} h(t,u)\left(\sqrt{F_n(u)-F_n(t)}-\sqrt{F_0(u)-F_0(t)}\right)^2 dt du\\
&= \int_{A_{1,n}} (h(t,t_0)+o(1))\frac{(u-t_0+\delta_n)^2f_0(t_0)^2+o(\delta_{n}^{2})}{4(F_0(t_0)-F_0(t))}dt du\\
&= \int_{0}^{t_0-\delta_n} h(t,t_0)\frac{f_0(t_0)^2(\delta_{n}^{3}+o(\delta_{n}^{3}))}{12(F_0(t_0)-F_0(t))}dt.
\end{align*}
The last integral can be split into two integrals over the sets
$[0,t_0-\kappa_n)$ and $[t_0-\kappa_n,t_0-\delta_n],$ where
$\kappa_n=(\log n)^{-1/3}.$ Since
$$
\int_{0}^{t_0-\kappa_n}
h(t,t_0)\frac{f_0(t_0)^2(\delta_{n}^{3}+o(\delta_{n}^{3}))}{12(F_0(t_0)-F_0(t))}dt=O(\delta_n^3\kappa_n^{-1})
$$
and
\begin{align*}
&\int_{t_0-\kappa_n}^{t_0-\delta_n}h(t,t_0)\frac{(\delta_{n}^{3}+o(\delta_{n}^{3}))f_0(t_0)^2}{12(F_0(t_0)-F_0(t))}dt\\
&= (f_0(t_0)(\delta_{n}^{3}+o(\delta_{n}^{3}))/12)\int_{t_0-\kappa_n}^{t_0-\delta_n} (h(t_0,t_0)+o(1))\frac{f_0(t)+o(1)}{(F_0(t_0)-F_0(t))}dt\\
&= (f_0(t_0)h(t_0,t_0)(\delta_{n}^{3}+o(\delta_{n}^{3}))/12)\left[-\log(F_0(t_0)-F_0(t))\right]_{t_0-\kappa_n}^{t_0-\delta_n}\\
&=f_0(t_0)h(t_0,t_0)c^3n^{-1}/36+o(n^{-1}),
\end{align*}
it follows that
\begin{align}\label{fin2}
&\int_{A_{1,n}} h(t,u)\left(\sqrt{F_n(u)-F_n(t)}-\sqrt{F_0(u)-F_0(t)}\right)^2 dt du\nonumber\\
&=f_0(t_0)h(t_0,t_0)c^3n^{-1}/36+o(n^{-1}).
\end{align}

Next, a straightforward computation shows that
\begin{align}
\label{fin3}
\int_{A_{1,n}}
&h(t,u)\left(\sqrt{1-F_n(u)}-\sqrt{1-F_0(u)}\right)^2)dt
du\nonumber\\
&=\int_{A_{1,n}}\frac{(u-t_0+\delta_n)^2f_0(t_0)^2}{4(1-F_0(t_0))}dt
du= O(\delta_n^3).
\end{align}
Using (\ref{fin}), (\ref{fin2}) and (\ref{fin3}), we get
\begin{align*}
&\int_{A_{1,n}}h(t,u)
\left(\sqrt{F_n(t)}-\sqrt{F_0(t)}\right)^2dtdu\\
&\qquad+\int_{A_{1,n}}
\left(\sqrt{F_n(u)-F_n(t)}-\sqrt{F_0(u)-F_0(t)}\right)^2dtdu\\
&\qquad+\int_{A_{1,n}}\left(\sqrt{1-F_n(u)}-\sqrt{1-F_0(u)}\right)^2dt
du\\
&=f_0(t_0)h(t_0,t_0)n^{-1}/36+O(\delta_n^3\kappa_n^{-1}).
\end{align*}
The integrals over $A_{2,n},A_{3,n}$ and $A_{4,n}$ can be treated
in a similar way. Indeed,
\begin{align*}
&\int_{A_{k,n}}h(t,u)
\left(\sqrt{F_n(t)}-\sqrt{F_0(t)}\right)^2dtdu\\
&\qquad+\int_{A_{k,n}}
\left(\sqrt{F_n(u)-F_n(t)}-\sqrt{F_0(u)-F_0(t)}\right)^2dtdu\\
&\qquad+\int_{A_{k,n}}\left(\sqrt{1-F_n(u)}-\sqrt{1-F_0(u)}\right)^2dt\,du\\
&=f_0(t_0)h(t_0,t_0)n^{-1}/36+O(\delta_n^3\kappa_n^{-1}),\quad
k=2,3,4.
\end{align*}

Moreover, it is easily verified that
\begin{align*}
&\int_{A_{5,n}}h(t,u)
\left(\sqrt{F_n(t)}-\sqrt{F_0(t)}\right)^2dtdu\\
&\qquad+\int_{A_{5,n}}
\left(\sqrt{F_n(u)-F_n(t)}-\sqrt{F_0(u)-F_0(t)}\right)^2dtdu\\
&\qquad
+\int_{A_{5,n}}\left(\sqrt{1-F_n(u)}-\sqrt{1-F_0(u)}\right)^2dt
du=O\big(\delta_n^3\big).
\end{align*}
Thus, we infer that the asymptotic squared Hellinger distance
between $q_0$ and $q_n$ is given by
$$
H^2(q_0,q_n)=f_0(t_0)h(t_0,t_0)n^{-1}/18.
$$
By using lemma \ref{mm} we now get:
\begin{align*}
&(n\log n)^{1/3} \inf_{T_n} \max \{E_{n,q_0}|T_n-F_0(t_0)|,E_{n,q_n}|T_n-F_n(t_0)|\}\\
&\ge \frac{1}{4} (n\log n)^{1/3}|F_n(t_0)-F_0(t_0)|\{1-H^2(q_n,q_0)\}^2\\
&\to\frac{1}{4}cf_0(t_0)\exp\left\{-\frac{1}{18}h(t_0,t_0)f(t_0)c^3
\right\}
\end{align*}
Maximizing the last expression over c yields the desired minimax
lower bound. \hspace*{0.5cm} \eop

\section{Asymptotic distribution of Birg\'e's estimator in the non-separated case}
\label{sec:Birge}
\setcounter{equation}{0}
\cite{lucien:99}  constructed a histogram-type estimator to show that the minimax lower bound rate of the preceding section can indeed be attained in the non-separated case. It is defined in the following way. Let $t_0$
be an interior point of [0,1], let $c$ be a positive constant and
let $K=\lfloor c^{-1}(n\log n)^{1/3}\rfloor,$ where $n$ is the
sample size and where $\lfloor x \rfloor$ denotes the ``floor" of
$x$, i.e., the largest integer which is smaller than or equal to
$x.$ We distinguish two cases.

\begin{enumerate}
\item[(i)]
If $Kt_0 \in \mathbb{N},$ the
interval $[0,1]$ is partitioned into $K$ intervals $I_j,$
$j=1,\dots,K,$ of equal length $1/K$, where $I_j=[t_j,t_{j+1})$,
$1\le j<K,$ $I_K=[t_K,t_{K+1}]$, and $t_1=0,\,t_{K+1}=1$.
\item[(ii)] If $Kt_0\notin \mathbb{N},$ the interval $[0,1]$ is partitioned into $K+1$
intervals $I_j$, where $I_j=[t_j,t_{j+1}),$ $1\le j\le K,$
$I_{K+1}=[t_{K+1},t_{K+2}]$, and $t_1=0,\,t_j=t_0-\left(\lfloor t_0K
\rfloor -j\right)/K$, $1\le j\le K+1$,
$t_{K+2}=1.$ Note that in this case the intervals $I_2,\ldots,I_K$
have length $1/K,$ but that $I_1$ and $I_{K+1}$ have a shorter
length. Furthermore, just as in case (i), $t_0$ is the left
boundary point of one of the intervals $I_j$.
\end{enumerate}

In fact we
slightly modified the definition of Birg\'e who always partitions the
interval into $K$ subintervals of equal length. The reason for our modification is that we
want to assign a fixed position to $t_0$ with respect to the
boundary points of the interval $I_j$ to which it belongs, since
the bias of the estimator heavily depends on this position.
Letting $t_0$ be a left boundary point enables us to compare the
results for different sample sizes ``on equal footing", so to
speak.

Let $\dd_{i,1}$, $\dd_{i,2}$ and $\dd_{i,3}$ be defined by (\ref{def_indicators}). We define, following \cite{lucien:99}, for
$1\le j,k\le K$,
\begin{align*}
&N_j=\#\left\{T_i:T_i\in I_j\right\},\quad M_j=\#\left\{U_i:U_i\in
I_j\right\}\\
&Q_{j,k}=\#\left\{(T_i,U_i):T_i\in I_j,\,U_i\in
I_k\right\},
\end{align*}
and
$$
N_j'=\sum_{T_i\in I_j}\dd_{i,1},\quad Q_{j,k}'=\sum_{T_i\in
I_j,\,U_i\in I_k}\dd_{i,2},\quad \quad M_j'=\sum_{U_i\in I_j}\dd_{i,3}\,.
$$
In addition to these (integer-valued) random variables,
\cite{lucien:99} defines the random variables:
\begin{equation}
\label{F_{j,k}}
\hat F^{(j,k)}=\left\{\begin{array}{lllll}
\displaystyle{\frac{N_k'}{N_k}-\frac{Q_{j,k}'}{Q_{j,k}}}&,\,j<k,\\
\displaystyle{1-\frac{M_k'}{M_k}+\frac{Q_{k,j}'}{Q_{k,j}}}&,\,j>k,
\end{array}
\right.
\end{equation}
weights $w_{j,k}$, defined by
\begin{equation}
\label{w_{j,k}}
w_{j,k}=\left\{\begin{array}{lllll}
\displaystyle{\frac{\sqrt{N_k\wedge(K Q_{j,k})}}{(k-j+1)W_j}}&,\,j<k,\\
\phantom{bla}&\\
\displaystyle{\frac{\sqrt{M_k\wedge(K Q_{k,j})}}{(j-k+1)W_j}}&,\,j>k,
\end{array}
\right.
\end{equation}
where
\begin{equation}
\label{def_W_j}
W_j=\sum_{k<j}\frac{\sqrt{M_k\wedge(KQ_{j,k})}}{j-k+1}+\sum_{k>j}\frac{\sqrt{N_k\wedge(KQ_{j,k})}}{k-j+1}\,.
\end{equation}

We are now ready to define Birg\'e's estimator $\widetilde{F}_n$.

\begin{definition}
{\rm
\label{def:Birge_est}
({\bf Birg\'e's estimator})
Let the intervals $I_j$ be defined as in (i) or (ii) above (depending on the value of
$t_0$), and let $\hat F^{(j,k)}$ and the weights $w_{j,k}$ be defined by (\ref{F_{j,k}}) and
(\ref{w_{j,k}}), respectively. Then, for $t$ belonging to the interval $I_j$, Birg\'e's estimator $\widetilde{F}_n(t)$ of $F_0(t)$ is
defined by
\begin{equation}
\label{Birge_est}
\widetilde{F}_n(t)=\sum_{k:k\ne j}w_{j,k}{\hat F^{(j,k)}}\,.
\end{equation}
}
\end{definition}

\vspace{0.3cm}
In determining the asymptotic distribution of Birg\'e's estimator, we are faced with the following difficulties.

\begin{enumerate}
\item[(1)]
The weights $w_{j,k}$ are ratios of random variables, which interact with the random variables
$M_k'/M_k$, $N_k'/N_k$ and $Q_{j,k}'/Q_{j,k}$, for which they are multipliers.
\item[(2)] The ratios $M_k'/M_k$, $N_k'/N_k$ and $Q_{j,k}'/Q_{j,k}$ are themselves
ratios of random variables.
\item[(3)] The weighted sum, defining Birg\'e's estimator, consists of {\it dependent}
summands. The dependence is caused by the dependence of the weights, the dependence
between the $M_k'/M_k$, $N_k'/N_k$ and $Q_{j,k}'/Q_{j,k}$ and the dependence between the weights and
these terms. This prevents a straightforward use of the
Lindeberg-Feller central limit theorem.
\end{enumerate}

These difficulties have to be dealt with in turn. The following crucial lemma bears on difficulty (1),
by showing that the random weights $w_{j,k}$ are close to deterministic weights $\widetilde{w}_{j,k}$.

\begin{lemma}
\label{lemma:asy}
Consider a partition of $[0,1]$ into $K$ or $K+1$ subintervals, according to the construction of Birg\'e's estimator,
using the scheme of (i) and (ii) at the beginning of this section. Assume that, for a fixed constant $c>0$,
\begin{equation}
\label{def_K_c}
K=K_n\sim \frac{(n\log n)^{1/3}}{c},\,n\to\infty,
\end{equation}
that is: the asymptotic binwidth is given by $c(n\log n)^{-1/3}$.
Moreover, assume that the observation density $h$ is continuous on the upper triangle of the unit square, staying away from
zero on its support. Let $g_1$ and $g_2$ be the first and second marginal density of $h$, respectively.
Finally, let $t_0$ be the left boundary point of $I_j$, let $a(t)$ and $b(t)$ be
defined by
\begin{equation}
\label{def_a_b}
a(t)=\sqrt{h(t_0,t)\wedge g_1(t)},\,\qquad
b(t)=\sqrt{h(t,t_0)\wedge g_2(t)},
\end{equation}
and let the deterministic weights $\widetilde{w}_{j,k}$ be defined by:
\begin{equation}
\label{tw_{j,k}}
\widetilde{w}_{j,k}=\left\{\begin{array}{lllll}
\displaystyle{\frac{3a(t_k)}{\left\{a(t_0)+b(t_0)\right\}(k-j+1)\log n}}&,k>j,\\
\phantom{bla}&\\
\displaystyle{\frac{3b(t_k)}{\left\{a(t_0)+b(t_0)\right\}(j-k+1)\log n}}&,\,k<j.
\end{array}
\right.
\end{equation}
Then:
\begin{enumerate}
\item[(i)]
\begin{equation}
\label{uniformity}
\sup_{k\ne j}\,(1+|j-k|)E\left|w_{j,k}-\widetilde{w}_{j,k}\right|=o\left(1/\log n\right),\,n\to\infty.
\end{equation}
\item[(ii)] $W_j$, defined by (\ref{def_W_j}), satisfies
\begin{equation}
\label{W_j-asymp}
W_j=\tfrac13(\log n)\sqrt{n/K}\left\{a(t_0)+b(t_0)\right\}\left\{1+o_p(1)\right\},\,n\to\infty,
\end{equation}
and, for $m=1,2,\dots$
\begin{equation}
\label{E1/W_j^m}
E\left\{1/W_j^m\right\}1_{\{W_j>0\}}\sim (9K/n)^{m/2}\left\{\left(a(t_0)+b(t_0)\right)\log n\right\}^{-m},\,n\to\infty.
\end{equation}
\end{enumerate}
\end{lemma}

\vspace*{0.5 cm}
It may be helpful to give some motivation for the construction of Birg\'e's statistic. If we replace $N_k,N_k'$, etc. by their expected
values, we obtain:
\begin{align*}
&\sum_{k>j}w_{j,k}\left\{\frac{\int_{I_k}F_0(u)\,dG_1(u)}{G_1(t_{k+1})-G_1(t_k)}-
\frac{\int_{t\in I_j,\,u\in I_k}\{F_0(u)-F_0(t)\}\,dH(t,u)}{\int_{t\in I_j,\,u\in I_k}\,dH(t,u)}\right\}\\
&+\sum_{k<j}w_{j,k}\left\{1-\frac{\int_{I_k}\{1-F_0(t)\}\,dG_2(u)}{G_2(t_{k+1})-G_2(t_k)}\right.\\+
&\left.\qquad\qquad+\frac{\int_{t\in I_k,\,u\in I_j}\{F_0(u)-F_0(t)\}\,dH(t,u)}{\int_{t\in I_k,\,u\in I_j}\,dH(t,u)}\right\},
\end{align*}
where $G_1$ and $G_2$ are the first and second marginal distribution functions of $H$, repectively. By expanding $F_0$ at the left endpoints
$t_k$ of the intervals, we get:
\begin{align}
\label{heuristic_expansion}
&\sum_{k>j}w_{j,k}\left\{F_0(t_k)-\{F_0(t_k)-F_0(t_j)\}\right\}\nonumber\\
&\qquad+\sum_{k<j}w_{j,k}\left\{1-\{1-F_0(t_k)\}+\{F_0(t_j)-F_0(t_k)\}\right\}
\nonumber\\
&\qquad+\frac1{2K^2}\sum_{k>j}w_{j,k}\left\{\frac{f_0(t_k)g_1(t_k)}{G_1(t_{k+1})-G_1(t_k)}-
\frac{\{f_0(t_k)-f_0(t_j)\}\,h(t_j,t_k)}{K\int_{t\in I_j,\,u\in I_k}\,dH(t,u)}\right\}\nonumber\\
&\qquad+\frac1{2K^2}\sum_{k<j}w_{j,k}\left\{\frac{f_0(t_k)g_1(t_k)}{G_2(t_{k+1})-G_2(t_k)}+
\frac{\{f_0(t_j)-f_0(t_k)\}\,h(t_k,t_j)}{K\int_{t\in I_k,\,u\in I_j}\,dH(t,u)}\right\}+\dots\nonumber\\
&=F_0(t_k)\sum_{k:k\ne j}w_{j,k}\nonumber\\
&\qquad+\frac1{2K}\sum_{k>j}w_{j,k}\left\{\frac{f_0(t_k)g_1(t_k)}{g_1(t_k)}-
\frac{\{f_0(t_k)-f_0(t_j)\}\,h(t_j,t_k)}{h(t_j,t_k)}\right\}\nonumber\\
&\qquad+\frac1{2K}\sum_{k<j}w_{j,k}\left\{\frac{f_0(t_k)g_1(t_k)}{g_2(t_k)}+
\frac{\{f_0(t_j)-f_0(t_k)\}\,h(t_k,t_j)}{h(t_k,t_j)}\right\}+\dots\nonumber\\
&=F_0(t_j)+\frac1{2K}\sum_{k>j}w_{j,k}\left\{f_0(t_k)-\{f_0(t_k)-f_0(t_j)\}\right\}\nonumber\\
&\qquad+\frac1{2K}\sum_{k<j}w_{j,k}\left\{f_0(t_k)+f_0(t_j)-f_0(t_k)\right\}+\dots\nonumber\\
&=F_0(t_j)+\frac1{2K}f_0(t_j)+\dots
\end{align}
One of the difficulties in this expansion that we have glossed over for the moment is that $g_1(t_k)$ tends to zero, if
$t_k\to1$, and that similarly $g_2(t_k)$ tends to zero, if $t_k\to0$. This difficulty has to be dealt with separately. We do not have that difficulty for $h$, since we assume that $h$ stays away from zero on its support.

The expansion suggests that the asymptotic bias at $t_j$ will be $f_0(t_j)/(2K)$, which is indeed the case. However, the expansion does not
explain the particular choice of the weights. Considering the deterministic counterparts $\widetilde{w}_{j,k}$ of $w_{j,k}$, given by
(\ref{tw_{j,k}}) in Lemma \ref{lemma:asy}, we see that the weights are proportional to $1/(1+|j-k|)$, which has the effect that the
smaller observation intervals give the biggest contribution to the estimator, taking advantage of the fact that the smaller observation
intervals do indeed give more precise information on the ``unobservable" $X_i$, if we know that $X_i$ is contained in the interval (see the
discussion on this point in section \ref{sec:intro}. The choice of these weights reduces the variance of the estimator. Only this fact is
responsible for the fact that the rate of convergence is slightly faster than $n^{-1/3}$.

It seems that the MLE is doing something similar automatically, but in a more efficient way, if we believe the ``working hypothesis",
discussed in section \ref{sec:intro}. Assuming the truth of this ``working hypothesis", the asymptotic variance of the MLE only involves
the local joint density $h$ of $(T_i,U_i)$ at $(t_0,t_0)$ and the density $f_0(t_0)$ of $X_i$ at $t_0$, whereas the variance of Birg\'e's
estimator also involves the marginal densities of $(T_i,U_i)$, which do not appear in the local minimax lower bound, derived in section
\ref{sec:minimax}.

Also note that the partition, needed in the construction of Birg\'e's estimator, is dependent on an a priori knowledge of whether we are in the separated or non-separated case; in the non-separated case binwidths of order $(n\log n)^{-1/3}$ are taken (otherwise the higher rate $(n\log n)^{-1/3}$ would not be attained), and in the separated case binwidths of order $n^{-1/3}$ (taking 
$(n\log n)^{-1/3}$ would let the variance dominate the bias, as the sample size tends to infinity). For the computation of the maximum likelihood estimator (MLE), discussed in section \ref{sec:MLE}, it is not necessary to use a priori knowledge on the observation distribution; the MLE, considered as a histogram adapts automatically to the separated or non-separated case and will choose generally smaller binwidth for the non-separated case. This is one of the major advantages of the MLE over Birg\'e's estimator, apart from being monotone with values restricted to $[0,1]$.

\vspace*{0.5 cm}
Using the notation of Lemma \ref{lemma:asy} we can now formulate the main result for Birg\'e's estimator.

\begin{theorem}
\label{limdis}
Let the observation density $h$ satisfy the same condition as in Lemma \ref{lemma:asy}, and let $F_0$ have a continuous derivative $f_0$ on
$(0,1)$, satisfying $f_0(t_0)>0$. Furthermore, let $I_j^{(n)}$ be a subinterval, belonging to the partition of
$[0,1]$ into $K$ intervals, corresponding to  the construction of Birg\'e's estimator for a sample of size $n$, where $K$ is defined by
(\ref{def_K_c})  in Lemma \ref{lemma:asy}. Finally, let $\a_n$ be defined by
\begin{equation}
\label{alpha_n}
\a_n=(n\log n)^{-1/3},
\end{equation}
and let $t_j^{(n)}$ be
the left boundary point of $I_j^{(n)}$, for which we assume that it converges to an interior point $t_0\in(0,1)$, as $n\to\infty$. Then:
\begin{enumerate}
\item[(i)]
\begin{equation}
\label{Birge_limit}
\a_n^{-1}
\left\{\widetilde{F}_n\left(t_j^{(n)}\right)-F_0\left(t_j^{(n)}\right)\right\}\stackrel{{\cal D}}{\longrightarrow}
N\left(\tfrac12cf_0(t_0),\s_0^2\right),\,n\to\infty.
\end{equation}
where the right-hand side of (\ref{Birge_limit})
denotes a normal random variable, with expectation $\tfrac12cf_0(t_0)$ and variance
\begin{equation}
\label{sigma_0}
\s_0^2=\frac{3f_0(t_0)\left\{a(t_0)^2+b(t_0)^2\right\}}{ch(t_0,t_0)\left\{a(t_0)+b(t_0)\right\}^2}\,,
\end{equation}
and where $c$, $a(t_0)$ and $b(t_0)$ are defined by (\ref{def_K_c}) and (\ref{def_a_b}).
\item[(ii)]
\begin{equation}
\label{bias_lim}
\lim_{n\to\infty}\a_n^{-1}E\left\{\widetilde{F}_n(t_j)-F_0(t_j)\right\}=\tfrac12cf_0(t_0),
\end{equation}
and
\begin{equation}
\label{var_lim}
\lim_{n\to\infty}\a_n^{-2}\mbox{\rm var}\left\{\widetilde{F}_n\left(t_j^{(n)}\right)\right\}=\s_0^2.
\end{equation}
\end{enumerate}
\end{theorem}

Note that Theorem \ref{limdis} implies that the optimal value of $c$ is given by
$$
\frac32\left(\frac{9f_0(t)^4}{2h(t,t)^2}\right)^{1/3}\left(\frac{a(t)^2+b(t)^2}{\{a(t)+b(t)\}^2}\right)^{2/3}.
$$
This value of the constant was used in the simulations, reported below.

\section{Birg\'e's estimator in the separated case}
\label{sec:separated}
\setcounter{equation}{0}
We consider the asymptotic behavior of Birg\'e's estimator in
the separated case. This is mainly meant for illustrative purposes and
to give background to the simulation study. We
therefore do not aim to prove results in the widest generality and confine
our discussion to the case where the density $h$ of the observed pairs
$(T_i,U_i)$ has as support the triangle with vertices $(0,\e)$, $(0,1)$ and $(1-\e,1)$
and stays away from zero on its support, which is the situation we consider in the
simulation study.
In this case the faster rate $(n\log n)^{-1/3}$ is unattainable,
and we know that Birg\'e's estimator (and also the MLE) can only achieve the rate $n^{-1/3}$.
We therefore assume $K$ to be of order
$n^{1/3}$ and set $K=\lfloor c^{-1}n^{1/3} \rfloor.$

As in section \ref{sec:Birge} we introduce deterministic weights $\widetilde{w}_{j,k}$ to replace the random weights $w_{j,k}$.
Recall that, by definition,
\begin{equation}
\label{wese}
w_{j,k}=\left\{\begin{array}{lllll}
\displaystyle{\frac{\sqrt{N_k\wedge(K Q_{j,k})}}{(k-j+1)W_j}}&,\,j<k,\\
\phantom{bla}&\\
\displaystyle{\frac{\sqrt{M_k\wedge(K
Q_{k,j})}}{(j-k+1)W_j}}&,\,j>k,
\end{array}
\right.
\end{equation}
and
$$
W_j=\sum_{1\le k < j} \frac{\sqrt{M_k\wedge (KQ_{j,k})}}{j-k+1}
+\sum_{j<k\le K} \frac{\sqrt{N_k\wedge (KQ_{j,k})}}{k-j+1}.
$$
Let $g_1$ and $g_2$ be the first and second marginal density of $h$, respectively, that is:
\begin{equation}
\label{marg_dens_sep}
g_1(t)=\int_t^1 h(t,u)\,du,\,g_2(t)=\int_0^t h(t',t)\,dt',\,t\in[0,1].
\end{equation}
Then, if $2\epsilon\le t_0 \le 1-2\epsilon$,
\begin{align*}
W_j&\sim \sum_{k:t_j-t_k>\epsilon}
\frac{\sqrt{cn^{2/3}\left\{h(t_k,t_j)\wedge g_2(t_k)\right\}}}{j-k+1}\\
&\qquad\qquad\qquad\qquad\qquad\qquad\qquad+ \sum_{k:t_k-t_j>\epsilon}
\frac{\sqrt{cn^{2/3}\left\{h(t_j,t_k)\wedge g_1(t_k)\right\}}}{k-j+1}\\
&\sim n^{1/3} \int_\epsilon^{t_0-\epsilon}
\frac{\sqrt{c\left\{h(t,t_0)\wedge g_2(t)\right\}}}{t_0-t}dt +n^{1/3} \int_{t_0+\epsilon}^{1-\epsilon}
\frac{\sqrt{c\left\{h(t_0,t)\wedge g_1(u)\right\}}}{t-t_0}dt,
\end{align*}
showing $W_j\asymp n^{1/3}$.
The deterministic weights $\widetilde{w}_{j,k}$ are now defined by:
\begin{equation}
\label{modified_w_sep}
\widetilde{w}_{j,k}=\left\{\begin{array}{l l}
\displaystyle{\frac{\sqrt{h(t_k,t_j)\wedge g_2(t_k)}}{K\widetilde{W}(t_0)\left(t_0-t_k\right)}\,,} & k<j,\\
\phantom{...} &\phantom{...}\\
\displaystyle{\frac{\sqrt{h(t_j,t_k)\wedge g_1(t_k)}}{K\widetilde{W}(t_0)\left(t_k-t_0\right)}\,,} &k>j,
\end{array}\right.
\end{equation}
where
\begin{equation}
\label{W_tilde_sep}
\widetilde{W}(t_0)=\int_{\e}^{t_0-\epsilon}
\frac{\sqrt{h(t,t_0)\wedge g_2(t)}}{t_0-t}\,dt +\int_{t_0+\epsilon}^{1-\e}
\frac{\sqrt{h(t_0,u)\wedge g_1(u)}}{u-t_0}\,du.
\end{equation}
We assume that the integrals on the right-hand side of (\ref{W_tilde_sep}) are finite, and hence that $\widetilde{W}(t_0)<\infty$.

We now have the following lemma, which plays a similar role as Lemma \ref{lemma:asy} in section \ref{sec:Birge}.

\begin{lemma}
\label{lemma:asy2}
Consider a partition of $[0,1]$ into $K$ or $K+1$ subintervals, according to the construction of Birg\'e's estimator,
using the scheme of (i) and (ii) at the beginning of section \ref{sec:Birge}. Assume that
$$
K=K_n\sim \frac{n^{1/3}}{c},\,n\to\infty,
$$
for a fixed constant $c>0$, that is: the asymptotic binwidth is given by $cn^{-1/3}$.
Let the weights $w_{j,k}$ and $\widetilde{w}_{j,k}$ be defined by (\ref{wese}) and (\ref{modified_w_sep}), respectively, where we assume
$\widetilde{W}(t_0)<\infty$. Then:
\begin{equation}
\label{uniformity_sep}
\sup_{k\ne j}\,(1+|j-k|)\left|w_{j,k}-\widetilde{w}_{j,k}\right|=o_p\left(n^{-1/3}\right),
\end{equation}
\end{lemma}

Using this lemma, we get the following limit result (compare with Theorem \ref{limdis}).

\begin{theorem}
\label{limdis_sep}
Suppose that the observation density $h$ has as support the triangle with vertices $(0,\e)$, $(0,1)$ and $(1-\e,1)$
and stays away from zero on its support. Let $F_0$ have a continuous derivative
$f_0$ on $(0,1)$, satisfying $f_0(t_0)>0$. Moreover, let $I_k^{(n)}$ be a subinterval, belonging to the partition of
$[0,1]$ into $K$ intervals, corresponding to  the construction of Birg\'e's estimator for a sample of size $n$.
Finally, let $\widetilde{W}(t_0)$ be defined by (\ref{W_tilde_sep}), where we assume  $\widetilde{W}(t_0)<\infty$.

Assume that, for a fixed constant $c>0$, $K=K_n\sim
n^{1/3}/c$, and let $t_k^{(n)}$ be the left boundary point of $I_k^{(n)}$, for which we assume that it converges to an interior point
$t_0\in(0,1)$, as
$n\to\infty$. Then we have, as $n\to \infty$
\begin{equation}
\label{Birge_limit_sep}
n^{1/3}
\left\{\widetilde{F}_n\left(t_k^{(n)}\right)-F_0\left(t_k^{(n)}\right)\right\}\stackrel{{\cal D}}{\longrightarrow}
N\left(\tfrac12cf_0(t_0),\s^2\right)
\end{equation}
where the right-hand side of (\ref{Birge_limit_sep})
denotes a normal random variable, with expectation $\tfrac12cf_0(t_0)$ and variance
\begin{align}
\label{as_var_sep}
&\s^2\nonumber\\
&=\frac{1}{c\widetilde{W}(t_0)^2}\int_{t_0+\epsilon}^{1-\epsilon}\frac{g_1(u)\wedge h(t_0,u)}{h(t_0,u)(u-t_0)^2}
\left\{F_0(u)-F_0(t_0)\right\}\left\{1-(F_0(u)-F_0(t_0))\right\}\,du\nonumber\\
&\quad+\frac{1}{c\widetilde{W}(t_0)^2}\int_{\epsilon}^{t_0-\epsilon}
\frac{g_2(t)\wedge h(t,t_0)}{h(t,t_0)(t_0-t)^2}\left\{F_0(t_0)-F_0(t)\right\}\left\{1-(F_0(t_0)-F_0(t))\right\}\,dt.
\end{align}
\end{theorem}

\vspace{0.3cm}
In the simulation study we take the observation density $h$ uniform on the triangle of its support. For ease of reference, we here
determine the value of the variance $\s^2$ of the asymptotic distribution for this case. If $h$ is uniform, its density is given by
\begin{equation}
\label{sep_obs_dens}
h(t,u)=\left\{\begin{array}{l l}\displaystyle{2(1-\e)^{-2},} & 0\le t+\epsilon\le u\le 1\\
0, & \textrm{elsewhere}
\end{array}\right..
\end{equation}
Hence the marginal densities $g_1$ and $g_2$ are given by:
$$
g_1(t)=\frac2{(1-\e)^2}\int_{t+\e}^1 du=\frac{2\{1-t-\e\}}{(1-\e)^2},\,t\in[0,1-\e],
$$
and
$$
g_2(u)=\frac2{(1-\e)^2}\int_0^{u-\e} du=\frac{2\{u-\e\}}{(1-\e)^2},\,u\in[\e,1].
$$
For $\widetilde{W}(t_0)$ we get:
\begin{align}
\label{W_tilde_sep_unif}
\widetilde{W}(t_0)&=\frac1{1-\e}\int_{\e}^{t_0-\epsilon}
\frac{\sqrt{2(t-\e)}}{t_0-t}\,dt +\frac1{1-\e}\int_{t_0+\epsilon}^{1-\e}
\frac{\sqrt{2(1-u-\e)}}{u-t_0}\,du.
\end{align}

Hence, using (\ref{as_var_sep}), we obtain:
\begin{align}
\label{sigma^2_uniform_sep}
\s^2&=\frac{1}{c\widetilde{W}(t_0)^2}\int_{t_0+\epsilon}^{1-\epsilon}\frac{1-u-\e}{(u-t_0)^2}
\left\{F_0(u)-F_0(t_0)\right\}\left\{1-(F_0(u)-F_0(t_0))\right\}\,du\nonumber\\
& \qquad
+\frac{1}{c\widetilde{W}(t_0)^2}\int_{\epsilon}^{t_0-\epsilon}
\frac{t-\e}{(t_0-t)^2}\left\{F_0(t_0)-F_0(t)\right\}\left\{1-(F_0(t_0)-F_0(t))\right\}\,dt.
\end{align}
where $\widetilde{W}(t_0)$ is defined by (\ref{W_tilde_sep_unif}).

\section{The maximum likelihood estimator}
\label{sec:MLE}
\setcounter{equation}{0}
As mentioned in section \ref{sec:intro}, the (nonparametric) maximum likelihood estimator (MLE or NPMLE) maximizes the (partial) log likelihood
$$
\sum_{i=1}^n \left\{\dd_{i1}\log F(T_i)+\dd_{i2}\log\left\{F(U_i)- F(T_i)\right\}+\dd_{i3}\log\left(1- F(U_i)\right)\right\},
$$
where the maximization is over all distribution functions $F$. For the non-sep\-arated case the following conjecture was given in \cite{Gr:91} (the lecture notes of a summer course given at Stanford University in 1990), which later appeared as part 2 of \cite{GrWe:92}):

\begin{theorem}
\label{th:nonsep_conjecture}
\mbox{\rm\bf (Conjecture in \cite{Gr:91})} Let $F_0$ and $H$ be continuously differentiable at $t_0$ and
$(t_0,t_0)$, respectively, with strictly positive derivatives $f_0(t_0)$ and $h(t_0,t_0)$, where $H$ is the distribution function of $(T_i,U_i)$. By continuous
differentiability of
$H$ at $(t_0,t_0)$ is meant that the density $h(t,u)$ is continuous at $(t,u)$, if $t<u$ and $(t,u)$ is
sufficiently close to $(t_0,t_0)$, and that $h(t,t)$, defined by
$$
h(t,t)=\lim_{u\downarrow t}h(t,u),
$$
is continuous at $t$, for $t$ in a neighborhood of $t_0$.

Let $0<F_0(t_0),H(t_0,t_0)<1$, and let $\hat F_n$ be the MLE of $F_0$. Then
$$
(n\log n)^{1/3}\left\{\hat F_n(t_0)-F_0(t_0)\right\}\bigm/\left\{\tfrac34f_0(t_0)^2/h(t_0,t_0)\right\}^{1/3}
\stackrel{{\cal D}}\longrightarrow 2Z,
$$
where $Z$ is the last time that standard two-sided Brownian motion minus the parabola $y(t)=t^2$ reaches its
maximum.
\end{theorem}

It was also shown in \cite{Gr:91} that Theorem \ref{th:nonsep_conjecture} is true for a ``toy" estimator, obtained by doing one step of the iterative convex minorant algorithm, starting the iterations at the underlying distribution function $F_0$; the ``toy" aspect is that we can of course not do this in practice. In spite of the fact that now more than 20 years have passed since this conjecture has been launched, it still has not been proved. In the simulation section we provide some material which seems to support the conjecture, but further research is necessary to settle this question.

For the separated case one can also introduce a toy estimator of the same type and one can again formulate the ``working hypothesis" that that the toy estimator and the MLE have the same pointwise limit behavior.
Anticipating that this would hold, \cite{We:95} derived the asymptotic distribution of the toy
estimator in the separated case, under the following conditions.

\begin{enumerate}
\item[(C1)] The support of $F_0$ is an interval $[0,M]$, where $M<\infty$.
\item[(C2)] $F_0$ and $H$ have densities $f_0$ and $h$ w.r.t.~Lebesgue measure on $\R$ and
$\R^2$, respectively.
\item[(C3)] Let the functions $k_{1,\e}$ and $k_{2,\e}$ be defined by
$$
k_{1,\e}(u)=\int_u^M\dfrac{h(u,v)}{F_0(v)-F_0(u)}\,\{F_0(v)-F_0(u)<\e^{-1}\}\,dv,
$$
and
$$
k_{2,\e}(v)=\int_0^v\dfrac{h(u,v)}{F_0(v)-F_0(u)}\,\{F_0(v)-F_0(u)<\e^{-1}\}\,du.
$$
Then, for $i=1,2$ and each $\e>0$,
$$
\lim_{\a\to\infty}\a\int_{(t_0,t_0+t/\a]}k_i(u,\e\a)\,du=0.
$$
\item[(C4)] $0<F_0(t_0)<1$ and $0<H(t_0,t_0)<1$.
\end{enumerate}
The motivation for these conditions is given in \cite{We:95} and actually become clear from the proof, which is not given here.

\begin{theorem}
\label{wellner}
\mbox{\rm (\cite{We:95})}
Suppose that assumptions (C1) to (C4) hold. Let $k_i,\,i=1,2,$ be
defined by
$$
k_1(u)=\int_u^M\dfrac{h(u,v)}{F_0(v)-F_0(u)}\,dv,\mbox{ and }
k_2(v)=\int_0^v\dfrac{h(u,v)}{F_0(v)-F_0(u)}\,du,
$$
and suppose that $f_0,g_1,g_2,k_1$ and $k_2$ are continuous at $t_0$, where $g_1$ and $g_2$ are the
first and second marginal densities of $h$, respectively. Moreover, assume
$f_0(t_0)>0$. Then, if $F_n^{(1)}$ is the estimator of the distribution function $F_0$, obtained
after one step of the iterative convex minorant algorithm, starting the iterations with $F_0$,
we have
$$
n^{1/3}\{2\xi(t_0)/f_0(t_0)\}^{1/3}\{F_n^{(1)}(t_0)-F_0(t_0)\}\stackrel{\D}\longrightarrow 2Z,
$$
where $Z$ is the last time where standard two-sided Brownian motion minus the parabola
$y(t)=t^2$ reaches its maximum, and where
$$
\xi(t_0)=\frac{g_1(t_0)}{F_0(t_0)}+k_1(t_0)+k_2(t_0)+\frac{g_2(t_0)}{1-F_0(t_0)}.
$$
\end{theorem}

It is indeed proved in \cite{piet:96} that, under slightly stronger conditions (the most important one being that an observation interval always has length $>\e$, for some $\e>0$), which hold for the examples in the simulation
below, the MLE has the same limit behavior, using the same norming constants. The expression for the
asymptotic variance in the separated case is remarkably different from the conjectured variance in the
non-separated case, which only depends on $F_0$ via
$f_0(t_0)$, showing that only the local behavior, depending on the density at $t_0$, is important
for the asymptotic variance (assuming that the working hypothesis holds).

Note that if $(T_i,U_i)$ is uniform on the upper triangle of the unit square, with vertices $(0,\e)$, $(0,1)$ and $(1-\e,1)$, we have:
$$
g_1(u)=\frac{2(1-u-\e)}{(1-\e)^2}\,,\qquad g_2(v)=\frac{2(v-\e)}{(1-\e)^2}\,,
$$
and, if $F_0$ is the uniform distribution function on $[0,1]$,
$$
k_1(u)=\frac{2\log\{(1-u)/\e\}}{(1-\e)^2}\,,\qquad k_2(v)=\frac{2\log(v/\e)}{(1-\e)^2}\,,
$$
so
$$
\xi(t_0)=\frac2{(1-\e)^2}\left\{\frac{1-t_0-\e}{t_0}+\log\left(\frac{t_0(1-t_0)}{\e^2}\right)+\frac{t_0-\e}{1-t_0}\right\}
$$
in this case. If $F_0$ is given by $F_0(x)=1-(1-x)^4,\,x\in[0,1]$, we get:
\begin{align}
\label{xi_for_sep}
&\xi(t_0)\\
&=\frac2{(1-\e)^2}\left\{\frac{1-t_0-\e}{F_0(t_0)}+\frac{t_0-\e}{1-F_0(t_0)}\right\}\nonumber\\
&\qquad+\frac1{2(1-\e)^2(1-t_0)^3}\left\{2\arctan\left(\frac{1-t_0-\e}{1-t_0}\right)+\log\left(\frac{2-2t_0-\e}{\e}\right)\right.\nonumber\\
&\qquad\left.+2\arctan\left(\frac{1-t_0+\e}{1-t_0}\right)-2\arctan\left(\frac1{1-t_0}\right)+\log\left(\frac{t_0\bigl(2-2t_0+\e\bigr)}{\e(2-t_0)}\right)\right\}.
\end{align}
We give some results for the latter model in section \ref{sec:simulations2}.

\section{The smoothed maximum likelihood estimator}
\label{sec:MSLE}
\setcounter{equation}{0}
Let $h$ be the density of $(T_i,U_i)$, with first marginal density $h_1$ and second marginal $h_2$, and let $\f_{t,b,F}$ be a solution of the integral equation (in $\f$):
\begin{align*}
\f(u)&=d_F(u)\left\{k_{t,b}(u)+\int_{v>u}\frac{\f(v)-\f(u)}{F(v)-F(u)}\,h(u,v)\,dv\right.\\
&\left.\qquad\qquad\qquad\qquad\qquad\qquad\qquad-
\int_{v<u} \frac{\f(u)-\f(v)}{F(u)-F(v)}\,h(v,u)\,dv\right\},
\end{align*}
where
$$
d_F(u)=\frac{F(u)\{1-F(u)\}}{h_1(u)\{1-F(u)\}+h_2(u)F(u)}\,,
$$
and the function $k_{t,b}$ is defined by
\begin{equation}
\label{k_t}
k_{t,b}(u)=b^{-1}K\left((t-u)/b\right).
\end{equation}
Moreover, let the function $\th_{t,b,F}$ be defined by
\begin{equation}
\label{def_theta}
\th_{t,b,F}(u,v,\d_1,\d_2)=-\frac{\d_1\f_{t,b,F}(u)}{F(u)}-\frac{\d_2\{\f_{t,b,F}(v)-\f_{t,b,F}(u)\}}{F(v)-F(u)}+
\frac{\d_3\f_{t,b,F}(v)}{1-F(v)}\,,
\end{equation}
where $u<v$. Then, as in \cite{GeGr:97} (separated case) and \cite{GeGr:99} (non-separated case), we have the representation
\begin{align*}
&\int \IK\left((t-u)/b\right)\,d\bigl(\hat F_n-F_0\bigr)(u)
=\int\th_{t,b,\hat F_n}(u,v,\d_1,\d_2)\,dP_0(u,v,\d_1\,d_2)\\
&=\int\frac{\f_{t,b,\hat F_n}(u)}{\hat F_n(u)}F_0(u)h_1(u)\,du\\
&\qquad\qquad+\int \frac{\f_{t,b,\hat F_n}(v)-\f_{t,b,\hat F_n}(u)}{\hat F_n(v)-\hat F_n(u)}\{F_0(v)-F_0(u)\}h(u,v)\,du\,dv\\
&\qquad\qquad-\int\frac{\f_{t,b,\hat F_n}(v)}{1-\hat F_n(v)}\{1-F_0(v)\}h_2(v)\,dv.
\end{align*}
For $F=F_0$ we get the integral equation:
\begin{align*}
\f(u)&=d_{F_0}(u)\left\{k_{t,b}(u)+\int_{v>u}\frac{\f(v)-\f(u)}{F_0(v)-F_0(u)}\,h(u,v)\,dv\right.\\
&\left.\qquad\qquad\qquad\qquad\qquad\qquad-\int_{v<u} \frac{\f(u)-\f(v)}{F_0(u)-F_0(v)}\,h(v,u)\,dv\right\}.
\end{align*}
Using the theory in \cite{GeGr:97} and \cite{GeGr:99} again, we get that the solution $\f_{t,b,F_0}$ gives as an approximation for $n\,\mbox{var}(\tilde F_n(t))$:
\begin{align*}
&E\,\th_{t,b,F_0}(T_1,U_1,\dd_{11},\dd_{12})^2\\
&=\int \frac{\f_{t,b,F_0}(u)^2}{F_0(u)}h_1(u)\,du\\
&\qquad+\int\frac{\bigl\{\f_{t,b,F_0}(v)-\f_{t,b,F_0}(u)\bigr\}^2}{F_0(v)-F_0(u)}h(u,v)\,du\,dv+\int \frac{\f_{t,b,F_0}(v)^2}{1-F_0(v)}h_2(v)\,dv.
\end{align*}
The approximation seems to work pretty well, as can be seen in table \ref{table:as_approx_smooth}, where we estimated the actual variance for samples of size $n=1000$ by generating $10,000$ samples of size $1000$ from a Uniform$(0,1)$ distribution $F_0$ and a uniform observation distribution $H$ on the upper triangle of the unit square.

\begin{figure}[!ht]
\begin{center}
\includegraphics[scale=0.45]{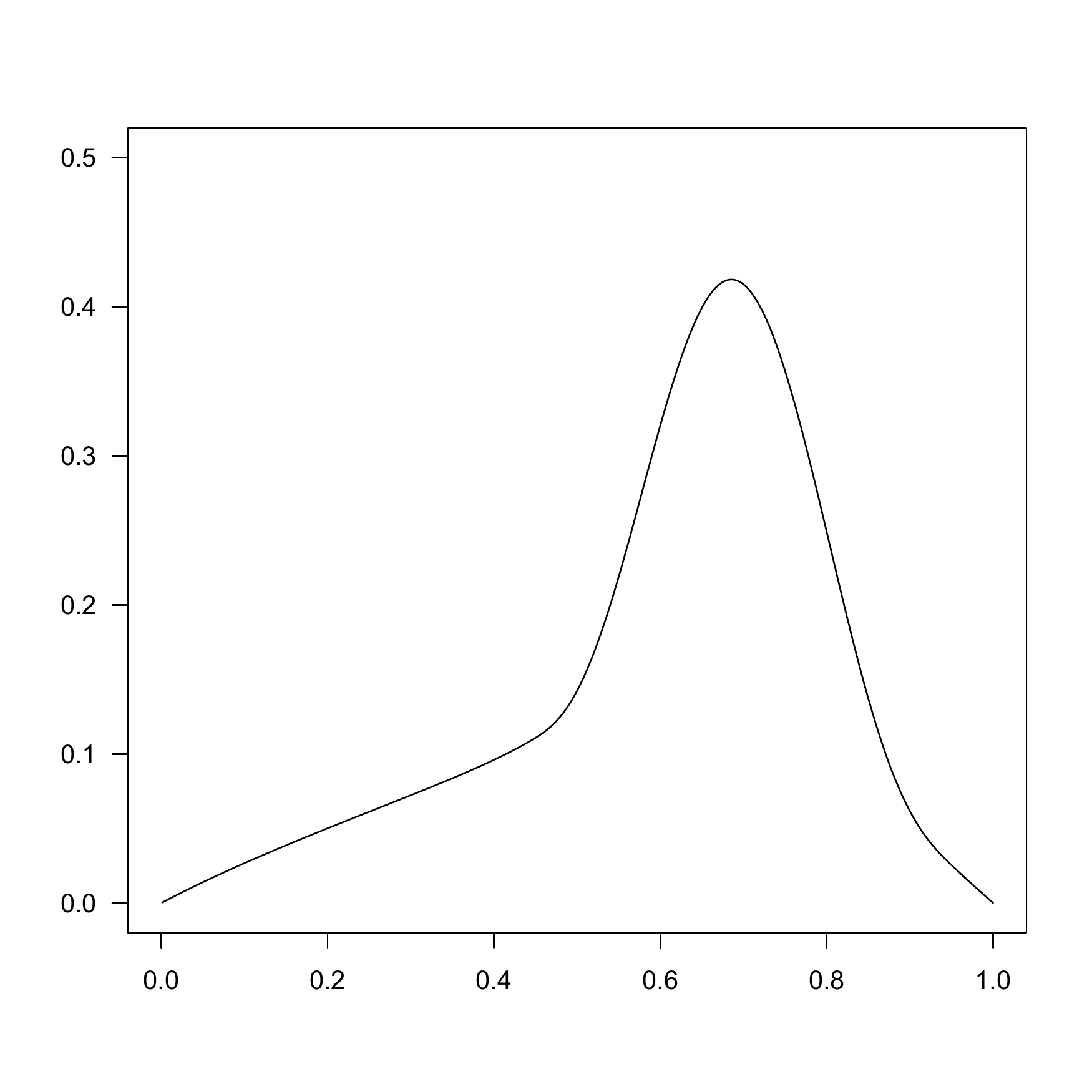}
\end{center}
\caption{The function $u\mapsto\f_{t,b_n,F_0}(u),\,u\in[0,1]$, for $t=0.7$, $b_n=n^{-1/5}$, $n=1000$, the Uniform$(0,1)$ distribution $F_0$ and a uniform observation distribution $H$ on the upper triangle of the unit square.}
\label{fig:phi}
\end{figure}

As in the papers cited above, we do not have an explicit expression for $\f_{t,b_n,F_0}$; a picture of $\f_{t,b_n,F_0}$ for $F_0$ the Uniform$(0,1)$ distribution $F_0$ and $b_n=n^{-1/5}$ is shown in Figure \ref{fig:phi}; the function was computed by solving the corresponding matrix equation on a $1000\times1000$ grid. Note that we apply the smooth functional theory of the above mentioned papers (which is also discussed in \cite{piet:96}) not for a fixed functional, but for changing functionals on shrinking intervals (in the hidden space). The reason we can do this is that the bandwidth $b$ is chosen to be of a larger order than the critical rate $n^{-1/3}$, and that then a different type of asymptotics sets in, with asymptotic normality, etc., instead of the non-standard asymptotics of the MLE itself. This method is also used in \cite{piet_geurt_birgit:10}, for the current status model.

\begin{table}
\caption{Estimates of the actual variances $\mbox{var}(\tilde F_n(t))$ (times $n$) and corresponding theoretical variances $E\th_{t,b_n,F_0}^2$, where $b_n=n^{-1/5}$, for sample size $n=1000$. The estimates of of the actual variances were based on $10,000$ samples of size $1000$ from a Uniform$(0,1)$ distribution $F_0$ and a uniform observation distribution $H$ on the upper triangle of the unit square.}
\label{table:as_approx_smooth}
\begin{tabular}{|l|c|c|c|c|c|c|c|c}
\hline
$t$  & $n\mbox{var}(\tilde F_n(t))$ &  $E\th_{t,b_n,F_0}^2$ & ratio\\
\hline
$0.1$ & 0.146489 & 0.142235 & 1.029910\\
$0.2$  &0.262056 & 0.255404 & 1.026044\\
$0.3$  &0.334990 & 0.332985 & 1.006019\\
$0.4$  &0.380357 & 0.376413 & 1.010479\\
$0.5$  &0.399258 & 0.390382 & 1.022736\\
$0.6$  &0.386292 & 0.376340 & 1.026444\\
$0.7$  &0.342651 & 0.332856 & 1.029428\\
$0.8$  &0.261457 & 0.255255 & 1.024296\\
$0.9$  &0.145304 & 0.142129 & 1.022338\\
\hline
\end{tabular}
\end{table}

In analogy with Theorem 4.2 in \cite{piet_geurt_birgit:10} we expect the following result to hold, using the conditions on the underlying distributions, discussed in \cite{GeGr:97} and \cite{GeGr:99}. To avoid messy notation, we will denote the smoothed MLE by $\tilde F_n$ instead of $\tilde F_n^{ML}$ in the remainder of this section.

\begin{theorem}
\label{th:limit_SMLE}
{\bf [Conjectured]}\,Let the conditions of Theorem 1, p.\ 212, in \cite{GeGr:97} (separated case) or Theorem 3.2, p.\ 647, in \cite{GeGr:99} (non-separated case) be satisfied. Moreover, let the joint density $h$ of the joint density of $(T_i,U_i)$ have a continuous bounded second total derivative in the interior of its domain and let $f_0$ have a continuous derivative at the interior point $t$ of the support of $f_0$, and let $\tilde F_n$ be the smoothed MLE, defined by (\ref{SMLE}). Then, if $b_n\asymp n^{-1/5}$, we have
$$
\sqrt{n}\left\{\tilde F_n(t)-F_0(t)-\tfrac12b_n^2f_0'(t)\int u^2K(u)\,du\right\}\Bigm/\s_n\stackrel{{\cal D}}\longrightarrow N\left(0,1\right),
\,n\to\infty,
$$
where $N(0,1)$ is the standard normal distribution and $\s_n^2$ is defined by
\begin{equation}
\label{var_SMLE}
\s_n^2= E\,\th_{t,b_n,F_0}\left(T_1,U_1,\dd_{11},\dd_{12}\right)^2,
\end{equation}
with $\th_{t,b_n,F_0}$ given by (\ref{def_theta}).
\end{theorem}

\begin{figure}[!ht]
\begin{center}
\includegraphics[scale=0.5]{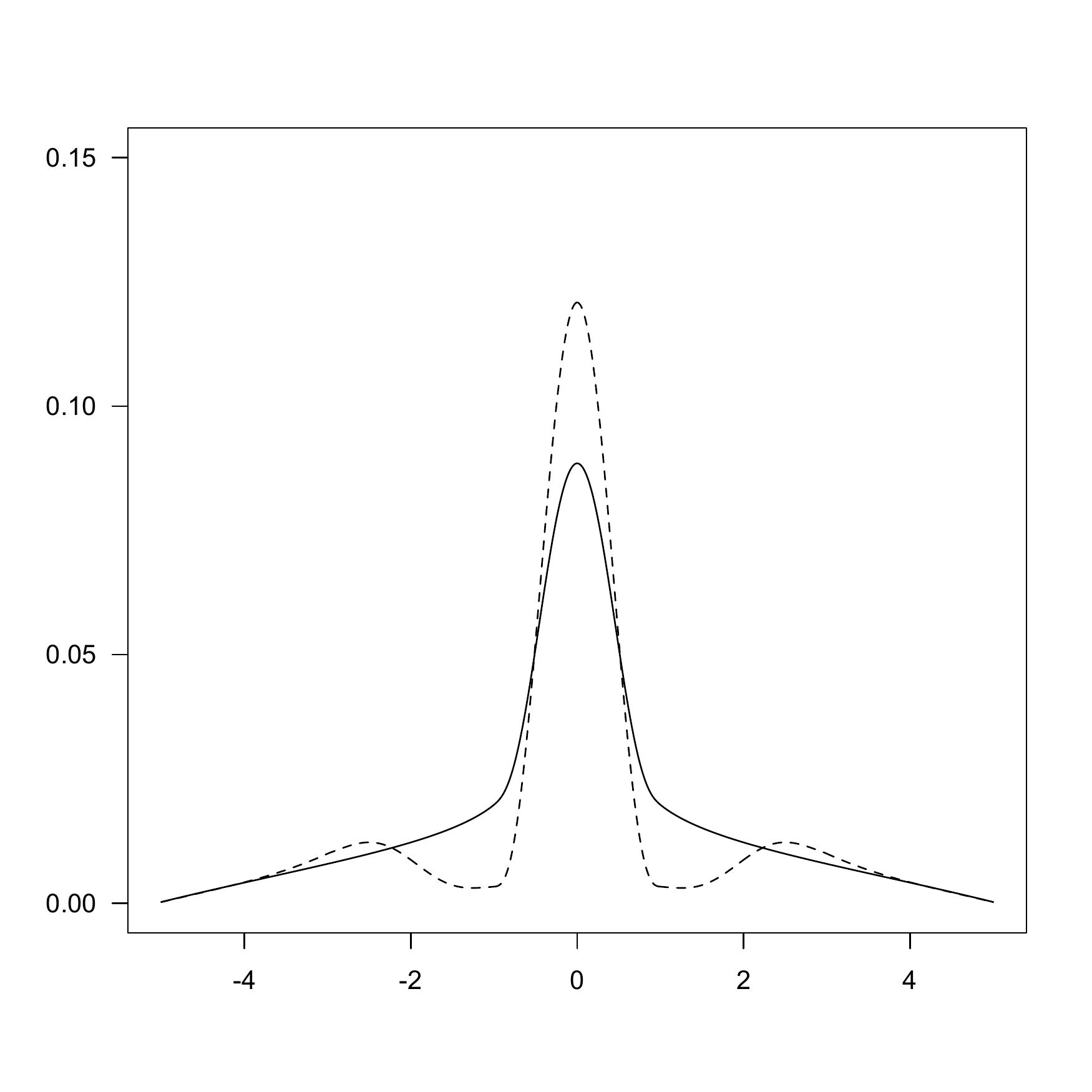}
\end{center}
\caption{The function $u\mapsto\f_{t,b,F_0}(t-b u),\,u\in[-5,5]$, for $t=0.5$, $b=0.1$, the Uniform$(0,1)$ distribution $F_0$ and (non-separated case:) a uniform observation distribution $H$ on the upper triangle of the unit square (solid curve) and the function $u\mapsto\f_{t,b,F_0}(t-b u)$ for the (separated) case where the observation distribution $H$ is uniform on the triangle with vertices $(0,\e)$, $(0,1)$ and $(1-\e,1)$, where $\e=0.2$ (dashed).}
\label{fig:two_phis}
\end{figure}

Note that (the conjectured) Theorem \ref{th:limit_SMLE} covers both the separated and the non-separated case. Unfortunately, we do not have an explicit expression for (\ref{var_SMLE}) in Theorem \ref{th:limit_SMLE} at present. The functions $\f_{F_0}$, defining the function $\th_{F_0}$ and hence also the variance $\s_n^2$, are of a rather different nature for the separated case and the non-separated case. For an example of this, see Figure \ref{fig:two_phis}.

The variance $\s_n^2$ can be estimated by
$$
\hat\s_n^2=\int \tilde\th_{t,b_n,\tilde F_n}(t,u,\d_1,\d_2)\,d\P_n(u,v,\d_1,\d_2),
$$
where
\begin{align*}
\tilde\th_{t,b_n,\tilde F_n}(u,v,\d_1,\d_2)
&=-\frac{\d_1\tilde\f_{t,b_n,\tilde F_n}(u)}{\tilde F_n(u)}-\frac{\d_2\{\tilde\f_{t,b_n,\tilde F_n}(v)-\tilde\f_{t,b_n,\tilde F_n}(u)\}}{\tilde F_n(v)-\tilde F_n(u)}\\
&\qquad\qquad\qquad\qquad\qquad+
\frac{\d_3\tilde\f_{t,b_n,\tilde F_n}(v)}{1-\tilde F_n(v)}\,,\,u<v,
\end{align*}
and $\tilde\f_{t,b_n,\tilde F_n}$ solves the integral equation
\begin{align}
\label{inteq2}
\f(u)&=d_{\tilde F_n(u)}(u)\left\{k_{t,b_n}(u)+\int_{v>u}\frac{\f(v)-\f(u)}{\tilde F_n(v)-\tilde F_n(u)}\,h_n(u,v)\,dv\right.\nonumber\\
&\left.\qquad\qquad\qquad\qquad\qquad\qquad-
\int_{v<u} \frac{\f(u)-\f(v)}{\tilde F_n(u)-\tilde F_n(v)}\,h_n(v,u)\,dv\right\},
\end{align}
and where $h_n$ is a kernel estimate of the density $h$, and where
\begin{align*}
&d_{\tilde F_n(u)}(u)=\frac{\tilde F_n(u)\{1-\tilde F_n(u)\}}{h_{n1}(u)\{1-\tilde F_n(u)\}+h_{n2}(u)\tilde F_n(u)},\\
&h_{n1}(u)=\int h_n(u,v)\,dv,\qquad h_{n2}(u)=\int h_n(v,u)\,dv.
\end{align*}
For $b_n$ chosen as in the theorem, the distribution function $\tilde F_n$ will be strictly increasing with probability tending to one. Since $\tilde F_n$ is also continuously differentiable, the equation (\ref{inteq2}) will have an absolutely continuous solution $\tilde\f_{t,b_n,\tilde F_n}$, and we do not have to take recourse to a solution {\it pair}, as in \cite{GeGr:99}, which deals separately with a discrete and absolutely continuous part.

In the corresponding result for the current status model we have explicit expressions, and we briefly discuss the analogy here, using a notation of the same type.
Let $\tilde F_n^{(CS)}$ be the smoothed MLE for the current status model, defined by (\ref{SMLE}), but now using the MLE $\hat F_n$ in the current status model. In this case the function $\th_{t,b,F}$, representing the functional in the observation space, is given by
\begin{equation}
\label{def_theta_CS}
\th_{t,b,F}^{(CS)}(u,\d)=-\frac{\d\f^{(CS)}_{t,b,F}(u)}{F(u)}+\frac{(1-\d)\f^{(CS)}_{t,b,F}(u)}{1-F(u)}\,,\,u\in(0,1).
\end{equation}
where $\f$ is given by:
$$
\f^{(CS)}_{t,b,F}(u)=\frac{F(u)\{1-F(u)\}}{g(u)}k_{t,b}(u),
$$
and $k_{t,b}$ is defined by (\ref{k_t}). Moreover, $g$ is the density of the (one-dimensional) observation distribution. The solution $\f_{t,b_n,F_0}^{(CS)}$ gives as an approximation for $n\,\mbox{var}(\tilde F_n(t))$:
\begin{align*}
&E\,\th^{(CS)}_{t,b_n,F_0}(T_1,\dd_1)^2=\int \frac{\f_{t,b_n,F_0}^{(CS)}(u)^2}{F_0(u)}g(u)\,du
+\int \frac{\f^{(CS)}_{t,b_n,F_0}(u)^2}{1-F_0(u)}g(u)\,du\\
&=\int \frac{F_0(u)\{1-F_0(u)\}k_{t,b_n}(u)^2}{g(u)}\,du\sim \frac{F_0(t)\{1-F_0(t)\}}{b_ng(t)}\int K(u)^2\,du,\,b_n\to0.
\end{align*}
Moreover,
\begin{align*}
\lim_{b\downarrow0}bE\,\th^{(CS)}_{t,b,F_0}(T_1,\dd_1)^2=\frac{F_0(t)\{1-F_0(t)\}}{g(t)}\int K(u)^2\,du,
\end{align*}
so in this case we obtain the central limit theorem
$$
\sqrt{n}\left\{\tilde F_n(t)-F_0(t)-\tfrac12b_n^2f_0'(t)\int u^2K(u)\,du\right\}\Bigm/\s_n\stackrel{{\cal D}}\longrightarrow N\left(0,1\right),
\,n\to\infty,
$$
where
$$
\s_n^2=E\,\th^{(CS)}_{t,b_n,F_0}(T_1,\dd_1)^2\sim\frac{F_0(t)\{1-F_0(t)\}}{b_ng(t)}\int K(u)^2\,du,
$$
see Theorem 4.2, p.\ 365, \cite{piet_geurt_birgit:10}.

\begin{remark}
{\rm It is tempting to think that the asymptotic variance can be found for case 2 by computing
$$
\lim_{b\downarrow0} bE\,\th_{t,b,F_0}\left(T_1,U_1,\dd_{11},\dd_{12}\right)^2,
$$
just as in the current status model. However, numerical computations suggested that $bE\,\th_{t,b,F_0}^2$ tends to zero in the non-separated case. This might mean that the variance is not of order $n^{-4/5}$ in this case, but perhaps contains a logarithmic factor, in analogy with the variance $(n\log n)^{-2/3}$ for the histogram-type estimators, like Birg\'e's estimator and the MLE without smoothing.

However, we do not expect this to happen for the separated case. All this still has to be determined by the analysis of the difference in asymptotic behavior of the functions $\f_{t,b_n,F_0}$ for the separated and non-separated case (see Figure \ref{fig:two_phis} for a picture of the rather different behavior of $\f_{t,b_n,F_0}$ in these two situations).
}
\end{remark}

\section{Simulation results for the non-separated case}
\label{sec:simulations}

In tables \ref{table2} to \ref{table_SMLE1} we present some simulation results for the
``non-separated case" for both Birg\'e's estimator, the MLE and the smoothed MLE. In
all cases the observation density was the uniform density on the
upper triangle. All results are based on 10,000 pseudo-random
samples. For Birg\'e's estimator the asymptotically optimal
binwidth was chosen in all simulations.

We study the case where $f_0$ is the uniform density on $[0,1]$
and give results for the interior points
$t_0=0.3,\,0.4,\,0.5$ and $0.6$. Although these points are
somewhat arbitrarily chosen, the results are representative for
what happens in the interior of the interval.

It can be seen from the tables that the squared bias for the MLE is,
in all cases, negligible compared to the variance. We note that
this is in contrast with Birg\'e's estimator. Moreover, the
variance of the MLE is generally smaller than that of Birg\'e's
estimator. Table \ref{table_SMLE1} shows, not unexpectedly, that the MSE of the smoothed MLE is much smaller than the MSE of either the MLE or Birg\'e's estimator.

\begin{table}
\caption{MSE for Birg\'e's estimator and MLE, times $(n\log n)^{2/3}$, $t_0=0.3,\,0.4,\,0.5$ and $0.6$, non-separated case.
The asymptotic MSE of Birg\'e's estimator and the conjectured MSE of the MLE are displayed in bold type.}
\label{table2}
\begin{center}
\begin{tabular}{l||c|c|c|c|c|c|c|c|}
 & \multicolumn{2}{|c|}{$t_0=0.3$} &
 \multicolumn{2}{|c|}{$t_0=0.4$}& \multicolumn{2}{c}{$t_0=0.5$} &
 \multicolumn{2}{|c}{$t_0=0.6$}\\
\cline{2-9} & Birg\'e & MLE & Birg\'e & MLE & Birg\'e & MLE &
Birg\'e & MLE\\
\cline{2-9}& {\bf 1.01} & {\bf 0.55} & {\bf 0.99} & {\bf 0.55} & {\bf 0.98} & {\bf 0.55} & {\bf 0.99} & {\bf 0.55}\\
\hline $n=1000$ & 1.10 & 0.50 & 1.10 & 0.55 & 1.09 & 0.55 & 1.11 & 0.55    \\
$n=2500$ & 1.06 & 0.52 & 1.08 & 0.54 & 1.07 & 0.55 & 1.06 & 0.53    \\
$n=5000$ & 1.05 & 0.50 & 1.03 & 0.54 & 1.04 & 0.56 & 1.03 & 0.53   \\
$n=10000$ & 1.03 & 0.51 & 1.02 & 0.54 & 1.00 & 0.54 & 1.06 & 0.54\\
\hline
\end{tabular}
\end{center}
\end{table}

\begin{table}
\caption{Variance for Birg\'e's estimator and MLE, times $(n\log n)^{2/3}$, $t_0=0.3,\,0.4,\,0.5$ and $0.6$, non-separated case. The asymptotic variance of Birg\'e's estimator and
the conjectured asymptotic variance of the MLE
 (MLE) are displayed in bold type.}
\label{table3}
\begin{center}
\begin{tabular}{l||c|c|c|c|c|c|c|c|}
 & \multicolumn{2}{|c|}{$t_0=0.3$} &
 \multicolumn{2}{|c|}{$t_0=0.4$}& \multicolumn{2}{c}{$t_0=0.5$} &
 \multicolumn{2}{|c}{$t_0=0.6$}\\
\cline{2-9} & Birg\'e & MLE & Birg\'e & MLE & Birg\'e & MLE &
Birg\'e & MLE\\
\cline{2-9}& {\bf 0.67} & {\bf 0.55} & {\bf 0.66} & {\bf 0.55} & {\bf 0.66} & {\bf 0.55} & {\bf 0.66} & {\bf 0.55}\\
\hline $n=1000$ & 0.79 & 0.50 & 0.78 & 0.55 & 0.78 & 0.55 & 0.79 & 0.55    \\
$n=2500$ & 0.75 & 0.52 & 0.75 & 0.54 & 0.74 & 0.55  & 0.73 & 0.53    \\
$n=5000$ & 0.74 & 0.50 & 0.71 & 0.54 & 0.73 & 0.56 & 0.72 & 0.53   \\
$n=10000$ & 0.69 & 0.51 & 0.69 & 0.54 & 0.69 & 0.54 & 0.72 & 0.54\\
\hline
\end{tabular}
\end{center}
\end{table}

\begin{table}
\caption{Squared Bias for Birg\'e's estimator and
MLE, times $(n\log n)^{2/3}$, $t_0=0.3,\,0.4,\,0.5$ and $0.6$, non-separated case. The asymptotic squared bias of Birg\'e's estimator is displayed in bold type.}
\label{table4}
\begin{center}
\begin{tabular}{l||c|c|c|c|c|c|c|c|}
 & \multicolumn{2}{|c|}{$t_0=0.3$} &
 \multicolumn{2}{|c|}{$t_0=0.4$}& \multicolumn{2}{c}{$t_0=0.5$} &
 \multicolumn{2}{|c}{$t_0=0.6$}\\
\cline{2-9} & Birg\'e & MLE & Birg\'e & MLE & Birg\'e & MLE &
Birg\'e & MLE\\
\cline{2-9}& {\bf 0.34} &  & {\bf 0.33} &  & {\bf 0.33} &  & {\bf 0.33} & \\
\hline $n=1000$ & 0.31 & $3.4\cdot 10^{-4}$ & 0.32 & $1.6\cdot 10^{-4}$ & 0.31 & $2.4\cdot 10^{-5}$ & 0.32 & $1.3\cdot 10^{-4}$    \\
$n=2500$ & 0.31 & $1.3\cdot 10^{-4}$ & 0.32 & $8.4\cdot 10^{-5}$ & 0.33 & $7.9\cdot 10^{-6}$  & 0.33 & $5.6\cdot 10^{-7}$    \\
$n=5000$ & 0.30 & $5.5\cdot 10^{-7}$ & 0.32 & $1.6\cdot 10^{-4}$ & 0.31 & $2.5\cdot 10^{-4}$ & 0.31 & $3.6\cdot 10^{-4}$   \\
$n=10000$ & 0.34 & $6.3\cdot 10^{-5}$ & 0.33 & $4.1\cdot 10^{-5}$
& 0.31 & $4.1\cdot 10^{-6}$ & 0.34 & $8.2\cdot 10^{-5}$\\
\hline
\end{tabular}
\end{center}
\end{table}

\begin{table}
\caption{MSE of SMLE divided by MSE of MLE, $t_0=0.3,\,0.4,\,0.5$ and $0.6$, non-separated case.}
\label{table_SMLE1}
\begin{center}
\begin{tabular}{l||c|c|c|c|}
 & \multicolumn{1}{|c|}{$t_0=0.3$} &
 \multicolumn{1}{c|}{$t_0=0.4$}& \multicolumn{1}{c|}{$t_0=0.5$} &
 \multicolumn{1}{c|}{$t_0=0.6$}\\
\cline{2-5}
 & ratio &  ratio &  ratio &  ratio\\
\hline $n=1000$  & 0.247  & 0.262 & 0.265  & 0.263\\
$n=2500$ & 0.217  & 0.236 & 0.236  & 0.233\\
$n=5000$ & 0.203  & 0.219 & 0.224  & 0.216\\
$n=10000$ & 0.187  & 0.197 & 0.204  & 0.201\\
\hline
\end{tabular}
\end{center}
\end{table}

\section{Simulation results for the separated case}
\label{sec:simulations2} 
For the separated case the results of a
simulation study are provided in the tables \ref{huge_samples} to \ref{table_SMLE3}. We first take again $F_0$ to be the uniform$(0,1)$ distribution function. On the other hand, we chose the
observation density defined by (\ref{sep_obs_dens}), with $\epsilon=0.1$, so the observation times $T_i$ and $U_i$ cannot become arbitrarily close.
The results are based on 10,000
pseudo-random samples. As in the non-separated case, the MSE of
the MLE turns out to be smaller than the MSE of Birg\'e's
estimator. Here the difference is however even more noticeable.

In the tables \ref{table6} to \ref{table8} we give the results for the MSE, variance and
squared bias for both estimators. Again it can be seen that the
variance of Birg\'e's estimator is generally larger than the
variance of the MLE. Moreover, as in the non-separated case,
the squared bias for the MLE is, in all cases, negligible
compared to the variance.

To show that the results are not specific for the uniform distribution, we give in the tables \ref{table10} to \ref{table12} the corresponding comparisons for the distribution function $F_0$, with density $f_0$, defined by
$$
F_0(x)=1-(1-x)^4,\qquad f_0(x)=4(1-x)^3,\qquad x\in[0,1].
$$
For the computation of the asymptotic variance of the MLE we used (\ref{xi_for_sep}) of section \ref{sec:MLE}. It is seen that the correspondence between the asymptotic expression for the variance and the actual sample variance of the MLE is rather good, and also that the superiority of the MLE w.r.t.\ Birg\'e's estimator is still more pronounced for this distribution function. Table \ref{table_SMLE3} shows that the ratio of the MSE of the SMLE and the MSE of the actual MLE is somewhat larger here, which is probably due to the fact that the asymptotic bias plays a larger role for the SMLE in this case (this bias vanishes for the uniform distribution function). The bias of the actual MLE is again very small for this distribution function, however.

As the fit with the asymptotic MSE was not satisfactory for Birg\'e's estimator in the separated case, we also did some simulations for much larger sample sizes. It turns out that the MSE then approximates the values predicted by the asymptotic theory. Some evidence is given in table \ref{huge_samples}. The
results are based on 1000 pseudo-random samples.

\begin{table}
\caption{MSE for Birg\'e's estimator divided by its asymptotic value, $t_0=0.3$, separated case.}
\label{huge_samples}
\begin{center}
\begin{tabular}{l||c|c|}
 & $f_0(t)=1$ & $f_0(t)=4(1-t)^3$\\
\cline{2-3}
$n=10^6$ & 1.12 & 1.09\\
$n=10^7$ & 1.04 & 1.04\\
\hline
\end{tabular}
\end{center}
\end{table}

\begin{table}
\caption{MSE for Birg\'e's estimator and MLE, times $n^{2/3}$,
$t_0=0.3,\,0.4,\,0.5$ and $0.6$, separated case. The asymptotic MSE (Birg\'e) and
``the asymptotic variance" (MLE) are displayed in bold type.}
\label{table6}
\begin{center}
\begin{tabular}{l||c|c|c|c|c|c|c|c|}
 & \multicolumn{2}{|c|}{$t_0=0.3$} &
 \multicolumn{2}{|c|}{$t_0=0.4$}& \multicolumn{2}{c}{$t_0=0.5$} &
 \multicolumn{2}{|c}{$t_0=0.6$}\\
\cline{2-9} & Birg\'e & MLE & Birg\'e & MLE & Birg\'e & MLE &
Birg\'e & MLE\\
\cline{2-9}& {\bf 0.34} & {\bf 0.12} & {\bf 0.32} & {\bf 0.13} & {\bf 0.31} & {\bf 0.13} & {\bf 0.32} & {\bf 0.13}\\
\hline $n=1000$ & 0.58 & 0.14 & 0.57 & 0.15 & 0.56 & 0.15 & 0.57 & 0.15    \\
$n=2500$ & 0.44 & 0.13 & 0.46 & 0.14 & 0.49 & 0.14 & 0.48 & 0.14    \\
$n=5000$ & 0.52 & 0.13 & 0.48 & 0.13 & 0.50 & 0.14 & 0.50 & 0.13   \\
$n=10000$ & 0.46 & 0.12 & 0.48 & 0.13 & 0.49 & 0.14 & 0.49 & 0.14\\
\hline
\end{tabular}
\end{center}
\end{table}

\begin{table}
\caption{Variance for Birg\'e's estimator and MLE, times $n^{2/3}$, $t_0=0.3,\,0.4,\,0.5$ and $0.6$, separated case. The asymptotic variances are
displayed in bold type.}
\label{table7}
\begin{center}
\begin{tabular}{l||c|c|c|c|c|c|c|c|}
 & \multicolumn{2}{|c|}{$t_0=0.3$} &
 \multicolumn{2}{|c|}{$t_0=0.4$}& \multicolumn{2}{c}{$t_0=0.5$} &
 \multicolumn{2}{|c}{$t_0=0.6$}\\
\cline{2-9} & Birg\'e & MLE & Birg\'e & MLE & Birg\'e & MLE &
Birg\'e & MLE\\
\cline{2-9}& {\bf 0.23} & {\bf 0.12} & {\bf 0.21} & {\bf 0.13} & {\bf 0.20} & {\bf 0.13} & {\bf 0.21} & {\bf 0.13}\\
\hline $n=1000$ & 0.46 & 0.14 & 0.47 & 0.15 & 0.46 & 0.15 & 0.47 & 0.15    \\
$n=2500$ & 0.35 & 0.13 & 0.36 & 0.14 & 0.39 & 0.14  & 0.37 & 0.14    \\
$n=5000$ & 0.42 & 0.13 & 0.38 & 0.13 & 0.40 & 0.14 & 0.39 & 0.13   \\
$n=10000$ & 0.36 & 0.12 & 0.37 & 0.14 & 0.39 & 0.14 & 0.39 & 0.14\\
\hline
\end{tabular}
\end{center}
\end{table}

\begin{table}
\caption{Squared Bias for Birg\'e's estimator and
MLE, times $n^{2/3}$, $t_0=0.3,\,0.4,\,0.5$ and $0.6$, separated case. The asymptotic squared bias (Birg\'e) is displayed in bold type.}
\label{table8}
\begin{center}
\begin{tabular}{l||c|c|c|c|c|c|c|c|}
 & \multicolumn{2}{|c|}{$t_0=0.3$} &
 \multicolumn{2}{|c|}{$t_0=0.4$}& \multicolumn{2}{c}{$t_0=0.5$} &
 \multicolumn{2}{|c}{$t_0=0.6$}\\
\cline{2-9} & Birg\'e & MLE & Birg\'e & MLE & Birg\'e & MLE &
Birg\'e & MLE\\
\cline{2-9}& {\bf 0.11} &  & {\bf 0.11} &  & {\bf 0.10} &  & {\bf 0.11} & \\
\hline $n=1000$ & 0.11 & $2.3\cdot 10^{-5}$ & 0.10 & $1.1\cdot 10^{-6}$ & 0.10 & $1.3\cdot 10^{-6}$ & 0.10 & $1.2\cdot 10^{-5}$    \\
$n=2500$ & 0.09 & $5.1\cdot 10^{-6}$ & 0.10 & $1.7\cdot 10^{-5}$ & 0.09 & $3.1\cdot 10^{-6}$  & 0.12 & $2.0\cdot 10^{-5}$    \\
$n=5000$ & 0.11 & $4.0\cdot 10^{-8}$ & 0.09 & $2.6\cdot 10^{-6}$ & 0.10 & $5.9\cdot 10^{-5}$ & 0.11 & $1.6\cdot 10^{-6}$   \\
$n=10000$ & 0.10 & $3.2\cdot 10^{-5}$ & 0.11 & $2.1\cdot 10^{-6}$
& 0.10 & $1.0\cdot 10^{-5}$ & 0.10 & $4.6\cdot 10^{-6}$\\
\hline
\end{tabular}
\end{center}
\end{table}

\begin{table}
\caption{MSE of SMLE divided by MSE of MLE,  $t_0=0.3,\,0.4,\,0.5$ and $0.6$,  separated case. }
\label{table_SMLE2}
\begin{center}
\begin{tabular}{l||c|c|c|c|}
 & \multicolumn{1}{|c|}{$t_0=0.3$} &
 \multicolumn{1}{c|}{$t_0=0.4$}& \multicolumn{1}{c|}{$t_0=0.5$} &
 \multicolumn{1}{c|}{$t_0=0.6$}\\
\cline{2-5}
 & ratio &  ratio &  ratio &  ratio\\
\hline $n=1000$  & 0.258  & 0.272 & 0.274  &  0.268\\
$n=2500$ & 0.230  & 0.244 & 0.243  & 0.244\\
$n=5000$ & 0.219  & 0.225 & 0.225  & 0.219\\
$n=10000$ & 0.199  & 0.201 & 0.206  & 0.203\\
\hline
\end{tabular}
\end{center}
\end{table}

\begin{table}
\caption{MSE for Birg\'e's estimator and MLE, times $n^{2/3}$, $f_0(t)=4(1-t)^3$, $t\in [0,1]$,
$t_0=0.3,\,0.4,\,0.5$ and $0.6$,  separated case. 
 The asymptotic MSE (Birg\'e) and
the asymptotic variance (MLE) are displayed in bold type.}
\label{table10}
\begin{center}
\begin{tabular}{l||c|c|c|c|c|c|c|c|}
 & \multicolumn{2}{|c|}{$t_0=0.3$} &
 \multicolumn{2}{|c|}{$t_0=0.4$}& \multicolumn{2}{c}{$t_0=0.5$} &
 \multicolumn{2}{|c}{$t_0=0.6$}\\
\cline{2-9} & Birg\'e & MLE & Birg\'e & MLE & Birg\'e & MLE &
Birg\'e & MLE\\
\cline{2-9}& {\bf 0.41} & {\bf 0.15} & {\bf 0.24} & {\bf 0.081} & {\bf 0.14} & {\bf 0.037} & {\bf 0.08} & {\bf 0.014}\\
\hline $n=1000$ & 0.53 & 0.16 & 0.39 & 0.088 & 0.21 & 0.041 & 0.101 & 0.016    \\
$n=2500$ & 0.61 & 0.16 & 0.33 & 0.087 & 0.25 & 0.039 & 0.100 & 0.015    \\
$n=5000$ & 0.56 & 0.16 & 0.36 & 0.083 & 0.18 & 0.038 & 0.101 & 0.014   \\
$n=10000$ & 0.49 & 0.15 & 0.36 & 0.082 & 0.22 & 0.037 & 0.120 & 0.014\\
\hline
\end{tabular}
\end{center}
\end{table}

\begin{table}
\caption{Variance for Birg\'e's estimator and MLE, times $n^{2/3}$, 
$t_0=0.3,\,0.4,\,0.5$ and $0.6,$
$f_0(t)=4(1-t)^3,\,t\in [0,1]$,  separated case.  The asymptotic
variances are displayed in bold type.}
\label{table11}
\begin{center}
\begin{tabular}{l||c|c|c|c|c|c|c|c|}
 & \multicolumn{2}{|c|}{$t_0=0.3$} &
 \multicolumn{2}{|c|}{$t_0=0.4$}& \multicolumn{2}{c}{$t_0=0.5$} &
 \multicolumn{2}{|c}{$t_0=0.6$}\\
\cline{2-9} & Birg\'e & MLE & Birg\'e & MLE & Birg\'e & MLE &
Birg\'e & MLE\\
\cline{2-9}& {\bf 0.28} & {\bf 0.15} & {\bf 0.16} & {\bf 0.081} & {\bf 0.091} & {\bf 0.037} & {\bf 0.051} & {\bf 0.014}\\
\hline $n=1000$ & 0.41 & 0.16 & 0.32 & 0.087 & 0.16 & 0.040 & 0.070 & 0.016    \\
$n=2500$ & 0.48 & 0.16 & 0.25 & 0.087 & 0.19 & 0.039  & 0.063 & 0.015    \\
$n=5000$ & 0.43 & 0.16 & 0.28 & 0.082 & 0.13 & 0.038 & 0.070 & 0.014   \\
$n=10000$ & 0.35 & 0.15 & 0.28 & 0.082 & 0.17 & 0.037 & 0.090 &
0.014\\
\hline
\end{tabular}
\end{center}
\end{table}

\begin{table}
\caption{Squared Bias for Birg\'e's estimator and
MLE, times $n^{2/3}$, $f_0(t)=4(1-t)^3$,\,$t\in [0,1]$, $t_0=0.3,\,0.4,\,0.5$ and $0.6$,  separated case. The asymptotic
squared bias (Birg\'e) is displayed in bold
type.}
\label{table12}
\begin{center}
\begin{tabular}{l||c|c|c|c|c|c|c|c|}
 & \multicolumn{2}{|c|}{$t_0=0.3$} &
 \multicolumn{2}{|c|}{$t_0=0.4$}& \multicolumn{2}{c}{$t_0=0.5$} &
 \multicolumn{2}{|c}{$t_0=0.6$}\\
\cline{2-9} & Birg\'e & MLE & Birg\'e & MLE & Birg\'e & MLE &
Birg\'e & MLE\\
\cline{2-9}& {\bf 0.14} &  & {\bf 0.079} &  & {\bf 0.045} &  & {\bf 0.025} & \\
\hline $n=1000$ & 0.12 & $1.1\cdot 10^{-4}$ & 0.076 & $1.6\cdot 10^{-4}$ & 0.051 & $2.2\cdot 10^{-4}$ & 0.030 & $3.2\cdot 10^{-4}$    \\
$n=2500$ & 0.13 & $3.2\cdot 10^{-5}$ & 0.080 & $2.3\cdot 10^{-4}$ & 0.054 & $2.0\cdot 10^{-4}$  & 0.037 & $1.1\cdot 10^{-4}$    \\
$n=5000$ & 0.13 & $1.4\cdot 10^{-6}$ & 0.075 & $3.0\cdot 10^{-4}$ & 0.048 & $1.0\cdot 10^{-4}$ & 0.031 & $8.7\cdot 10^{-5}$   \\
$n=10000$ & 0.13 & $4.8\cdot 10^{-5}$ & 0.079 & $1.4\cdot 10^{-4}$
& 0.049 & $1.1\cdot 10^{-4}$ & 0.030 & $8.0\cdot 10^{-5}$\\
\hline
\end{tabular}
\end{center}
\end{table}

\begin{table}
\caption{MSE of SMLE divided by MSE of MLE, times $n^{2/3}$, $f_0(t)=4(1-t)^3$,\,$t\in [0,1]$, $t_0=0.3,\,0.4,\,0.5$ and $0.6$, separated case,}
\label{table_SMLE3}
\begin{center}
\begin{tabular}{l||c|c|c|c|}
 & \multicolumn{1}{|c|}{$t_0=0.3$} &
 \multicolumn{1}{c|}{$t_0=0.4$}& \multicolumn{1}{c|}{$t_0=0.5$} &
 \multicolumn{1}{c|}{$t_0=0.6$}\\
\cline{2-5}
 & ratio &  ratio &  ratio &  ratio\\
\hline $n=1000$  & 0.439  & 0.395 & 0.443  &  0.435\\
$n=2500$ & 0.372  & 0.393 & 0.409  & 0.424\\
$n=5000$ & 0.350  & 0.354 & 0.383  & 0.391\\
$n=10000$ & 0.312  & 0.332 & 0.349  & 0.389\\
\hline
\end{tabular}
\end{center}
\end{table}

\section{Summary}
\label{sec:conclusion}
In the preceding, the limit distributions of three estimators for the interval censoring, case 2, problem were discussed: Birg\'e's estimator, the (nonparametric) maximum likelihood estimator (MLE) and the smoothed MLE, which is analogous to the smoothed MLE introduced in \cite{piet_geurt_birgit:10} for the current status model. Birg\'e's estimator is mainly of theoretical interest and constructed to show that the minimax rate can be attained. The construction uses prior knowledge on whether the observation distribution has arbitrarily small observation intervals (the so-called non-separated case) or not (the separated case). Such prior knowledge is not necessary for the MLE, which adapts automatically to either situation.

The conjectured limit distribution of the MLE in the non-separated case, given in \cite{Gr:91}, was (partially) checked in a simulation study, comparing Birg\'e's estimator, the MLE and the smoothed MLE. The simulation study seems to support the conjecture. The smoothed MLE converges at a faster rate than either Birg\'e's estimator or the MLE on which it is based if the underlying distribution is smooth, as is also borne out by the simulation study.

The limit distribution of the MLE in the separated case was given in \cite{piet:96} and the simulation study for the separated case shows that the asymptotic variance, arising from this result, provides a good approximation to the actual finite sample variance. The difference in behavior for the separated and non-separated cases persists for the smoothed MLE and in that case crucially depends on properties of the solution of an integral equation, as discussed in section \ref{sec:MSLE}. This analysis is based on a local version of the theory developed in \cite{GeGr:96}, \cite{GeGr:97} and \cite{GeGr:99}. The (numerical) solution of the integral equation can be used to estimate the variance of the smoothed MLE. The theoretically computed asymptotic variance, using a numerical solution of the integral equation, fits the observed sample variance rather well, but the discussion on this matter is heuristic and still contains lots of open questions.

\section{Appendix}

We split the proof of Theorem \ref{limdis} into several parts, dealing with the difficulties (1), (2) and (3),
mentioned in section \ref{sec:Birge}.
Here and in the following we will use some empirical process notation to make the
transition to the asymptotic distribution more transparant.
As an example, we give a representation of
\begin{equation}
\label{rep1}
N_k'/N_k=\frac{\sum_{T_i\in I_k}\dd_{i,1}}{\#\left\{T_i\in I_k\right\}}\,.
\end{equation}
in terms of integrals with respect to empirical distributions.
First we write:
$$
n^{-1}N_k'=\int_{t\in I_k}\d_1\,d\P_n(t,u,\d),
$$
where $\d=(\d_1,\d_2,\d_3)$ is the vector of indicators
$$
\d_1=1_{\{x\le t\}},\,\d_2=1_{\{t<x\le u\}},\,\d_3=1_{\{x>u\}},
$$
giving the position of the unobservable random variables $X_i$ with respect to the observation interval $[T_i,U_i]$,
and where $\P_n$ is the empirical measure of the random variables $(T_i,U_i,\dd_i)$
$=$ $(T_i,U_i,\dd_{i,1},\dd_{i,2},\dd_{i,3})$.

The denominator of (\ref{rep1}), after dividing by $n$, is rewritten in the form:
\begin{equation}
\label{N_k}
n^{-1}N_k=\int_{t\in I_k}\,d\G_{n,1}(t)=\G_{n,1}\left(t_{k+1}\right)-\G_{n,1}\left(t_k\right),
\end{equation}
where $\G_{n,1}$ is the empirical distribution function of the $T_i$, with underlying df $G_1$ and underlying density $g_1$, which is
the first marginal of $h$.
Using this notation, we get:
\begin{equation}
\label{rep1a}
N_k'/N_k=\frac{\int_{t\in I_k}\d_1\,d\P_n(t,u,\d)}{\G_{n,1}\left(t_{k+1}\right)-\G_{n,1}\left(t_k\right)}\,,
\end{equation}
where we the define the ratio to be zero if the denominator is zero. The terms $M_k'/M_k$ and $Q'_{j,k}/Q_{j,k}$ can be rewritten in a
similar way.

We will also use the following decomposition:
\begin{align}
\label{condit_exp_decomp1}
&\left\{\frac{N_k'}{N_k}-F_0(t_k)\right\}1_{\{N_k>0\}}\nonumber\\
&=\frac{N_k'-E\left\{N_k'|N_k\right\}}{N_k}\,1_{\{N_k>0\}}
+\frac{E\left\{N_k'-N_kF_0(t_k)|N_k\right\}}{N_k}\,1_{\{N_k>0\}}\,.
\end{align}
We similarly have, denoting $1-F_0$ by $\overline{F}_0$,
\begin{align}
\label{condit_exp_decomp2}
&\left\{\frac{M_k'}{M_k}-\overline{F}_0(t_k)\right\}1_{\{M_k>0\}}\nonumber\\
&=\frac{M_k'-E\left\{M_k'|M_k\right\}}{M_k}\,1_{\{M_k>0\}}
+\frac{E\left\{M_k'-M_k\overline{F}_0(t_k)|M_k\right\}}{M_k}\,1_{\{M_k>0\}}\,.
\end{align}
and
\begin{align}
\label{condit_exp_decomp3}
&\left\{\frac{Q_{j,k}'}{Q_{j,k}}-\left\{F_0(t_k)-F_0(t_j)\right\}\right\}1_{\{Q_{j,k}>0\}}\nonumber\\
&=\frac{Q_{j,k}'-E\left\{Q_{j,k}'|Q_{j,k}\right\}}{Q_{j,k}}\,1_{\{Q_{j,k}>0\}}\nonumber\\
&\qquad+\frac{E\left\{Q_{j,k}'-Q_{j,k}\left\{F_0(t_k)-F_0(t_j)\right\}|Q_{j,k}\right\}}{Q_{j,k}}\,1_{\{Q_{j,k}>0\}}\,.
\end{align}
One can consider this as a decomposition into a ``variance part" and a ``bias part", where the first terms on the right-hand sides of the
above expressions correspond to the variance part and the second terms to the bias part.

We first deal with the bias part.

\begin{lemma}
\label{lemma:bias}
Let the conditions of Theorem \ref{limdis} be satisfied, and let, for each interval $I_k$ of the partition, $t_k=t_k^{(n)}$ be its left
boundary point. Moreover, let $t_j=t_j^{(n)}\to t_0$, and $\a_n$ be defined by (\ref{alpha_n}).
Then we have for Birg\'e's statistic, defined by (\ref{Birge_est}),
\begin{enumerate}
\item[(i)] As $n\to\infty$,
\begin{equation}
\label{cond_var_lim}
\a_n^{-2}\,\mbox{\rm var}
\left(E\left\{\widetilde{F}_n(t_j)-F_0(t_j)\bigm|N_k,Q_{j,k},\,k>j;\,M_k,Q_{k,j},\,k<j\right\}\right)^2\to0.
\end{equation}
\item[(ii)] As $n\to\infty$,
\begin{equation}
\label{second_exp_lim}
\a_n^{-1}E\left\{\widetilde{F}_n(t_j)-F_0(t_j)\bigm|N_k,Q_{j,k},\,k>j;\,M_k,Q_{k,j},\,k<j\right\}\stackrel{p}\longrightarrow
\tfrac12cf_0(t_0).
\end{equation}
\end{enumerate}
\end{lemma}

\noindent
{\bf Proof.}\\
(i). If $N_k,\,M_k,\,Q_{j,k}$ and $Q_{k,j}$ are strictly positive, for all (relevant) values of $k$, we can write
\begin{align*}
&E\left\{\widetilde{F}_n(t_j)-F_0(t_j)\bigm|N_k,Q_{j,k},\,k>j;\,M_k,Q_{k,j},\,k<j\right\}\\
&=\sum_{k:k>j}w_{j,k}\left\{\frac{E\left\{N_k'-N_kF_0(t_k)|N_k\right\}}{N_k}-\right.\\
&\left.\qquad\qquad\qquad\qquad\qquad\qquad\frac{E\left\{Q_{j,k}'-Q_{j,k}\left\{F_0(t_k)-F_0(t_j)\right\}|Q_{j,k}\right\}}{Q_{j,k}}\right\}\\
&+\sum_{k:k<j}w_{j,k}\left\{\frac{E\left\{M_k-M_k'\overline{F}_0(t_k)|M_k\right\}}{M_k}\right.\\
&\left.\qquad\qquad\qquad\qquad\qquad\qquad+\frac{E\left\{Q_{k,j}'-Q_{k,j}\left\{F_0(t_k)-F_0(t_j)\right\}|Q_{k,j}\right\}}{Q_{k,j}}\right\},
\end{align*}
see (\ref{condit_exp_decomp1}) to (\ref{condit_exp_decomp3}). We can write this in the following form:
\begin{align*}
&E\left\{\widetilde{F}_n(t_j)-F_0(t_j)\bigm|N_k,Q_{j,k},\,k>j;\,M_k,Q_{k,j},\,k<j\right\}\\
&=\sum_{k:k>j}w_{j,k}\left\{\frac{\int_{I_k}\left\{F_0(t)-F_0(t_k)\right\}\,d\G_{n,1}(t)}{\G_{n,1}(t_{k+1})-\G_{n,1}(t_k)}\right.\\
&\left.\qquad\qquad\qquad\qquad-\frac{\int_{t\in I_j,\,u\in I_k}\left\{F_0(u)-F_0(t)-F_0(t_k)+F_0(t_j)\right\}\,d\H_n(t,u)}
{\int_{t\in I_j,\,u\in I_k}\,d\H_n(t,u)}\right\}\\
&\qquad+\sum_{k:k<j}w_{j,k}\left\{\frac{\int_{I_k}\left\{F_0(t)-F_0(t_k)\right\}\,d\G_{n,2}(t)}{\G_{n,2}(t_{k+1})-\G_{n,2}(t_k)}\right.\\
&\left.\qquad\qquad\qquad\qquad+\frac{\int_{t\in I_k,\,u\in I_j}\left\{F_0(u)-F_0(t)-F_0(t_j)+F_0(t_k)\right\}\,d\H_n(t,u)}
{\int_{t\in I_k,\,u\in I_j}\,d\H_n(t,u)}\right\}.
\end{align*}
By expanding $F_0$ in $t_k$ and $t_j$, as in (\ref{heuristic_expansion}), we find that this can be written
\begin{align*}
&E\left\{\widetilde{F}_n(t_j)-F_0(t_j)\bigm|N_k,Q_{j,k},\,k>j;\,M_k,Q_{k,j},\,k<j\right\}\\
&=\sum_{k:k>j}w_{j,k}\left\{\frac{f_0(t_k)\int_{I_k}\left(t-t_k\right)\,d\G_{n,1}(t)}{\G_{n,1}(t_{k+1})-\G_{n,1}(t_k)}\right.\\
&\left.\qquad\qquad\qquad\qquad-\frac{\int_{t\in I_j,\,u\in
I_k}\left\{f_0(t_k)\left(u-t_k\right)-f_0(t_j)\left(t-t_j\right)\right\}\,d\H_n(t,u)} {\int_{t\in I_j,\,u\in I_k}\,d\H_n(t,u)}\right\}\\
&\qquad+\sum_{k:k<j}w_{j,k}\left\{\frac{f_0(t_k)\int_{I_k}\left(t-t_k\right)\,d\G_{n,2}(t)}{\G_{n,2}(t_{k+1})-\G_{n,2}(t_k)}\right.\\
&\left.\qquad\qquad\qquad\qquad-\frac{\int_{t\in I_k,\,u\in
I_j}\left\{f_0(t_j)\left(u-t_j\right)-f_0(t_k)\left(t-t_k\right)\right\}\,d\H_n(t,u)} {\int_{t\in I_k,\,u\in I_j}\,d\H_n(t,u)}\right\}\\
&\qquad+o\left(1/K\right).
\end{align*}
The remainder term $o(1/K)$ arises from the fact that we can write, for example,
\begin{align*}
\frac{\int_{I_k}\left\{F_0(t)-F_0(t_k)\right\}\,d\G_{n,1}(t)}{\G_{n,1}(t_{k+1})-\G_{n,1}(t_k)}
=\frac{f_0(t_k)\int_{I_k}\left(t-t_k\right)\,d\G_{n,1}(t)}{\G_{n,1}(t_{k+1})-\G_{n,1}(t_k)}
+\left(t_{k+1}-t_k\right)o(1),
\end{align*}
where $t_{k+1}-t_k\le 1/K$, and the $o(1)$-factor is uniform in $k$, by the uniform contiuity of $f_0$ on $[0,1]$.
A similar expansion is used for the other terms, and the $o(1/K)$ remainder term now surfaces from the fact that the weights $w_{j,k}$ sum
to 1.

Furthermore, if $j<k$, and $t_k,t_j\in[\e,1-\e]$, for some $\e\in(0,1/2)$, we get:
\begin{align*}
&Ew_{j,k}^2\left\{\frac{f_0(t_k)\int_{I_k}\left(t-t_k\right)\,d\G_{n,1}(t)}{\G_{n,1}(t_{k+1})-\G_{n,1}(t_k)}
-\frac{f_0(t_k)\int_{I_k}\left(t-t_k\right)\,dG_1(t)}{\G_{n,1}(t_{k+1})-\G_{n,1}(t_k)}\right\}^21_{\{N_k>0>0\}}\\
&\le
E\frac{nf_0(t_k)^2\left\{\int_{I_k}\left(t-t_k\right)\,d\left(\G_{n,1}-G_1\right)(t)\right\}^2}
{(1+k-j)^2W_j^2\left\{\G_{n,1}(t_{k+1})-\G_{n,1}(t_k)\right\}}
1_{\{N_k>0\}}\\
&\sim \frac{f_0(t_k)^2}{3(1+k-j)^2Kg_1(t_k)}E\frac1{W_j^2}\,1{\{W_j>0\}}\\
&\sim\frac{3f_0(t_k)^2}{n(1+k-j)^2g_1(t_k)\left\{a(t_0)+b(t_0)\right\}^2(\log n)^2}\,,
\end{align*}
where we use (\ref{E1/W_j^m}) and exponential inequalities of the type discussed in the proof of Lemma \ref{lemma:asy} below for the probability that
$$
|\G_{n,1}(t_{k+1})-\G_{n,1}(t_k)-G_1(t_{k+1})+G_1(t_k)|>\e.
$$

We similarly get, for all $k>j$,
\begin{align*}
&E w_{j,k}^2\left\{\frac{\int_{t\in I_j,\,u\in
I_k}f_0(t_k)\left(u-t_k\right)\,d\left(\H_n-H\right)(t,u)}{\int_{t'\in I_j,\,u'\in I_k}\,d\H_n(t',u')}\right.\\
&\left.\qquad\qquad\qquad-\frac{\int_{t\in I_j,\,u\in
I_k}f_0(t_j)\left(t-t_j\right)\,d\left(\H_n-H\right)(t,u)}{\int_{t'\in I_j,\,u'\in I_k}\,d\H_n(t',u')}
\right\}^21_{\{Q_{j,k}>0\}}\\
&\le
E\frac{\left\{f_0(t_k)^2+f_0(t_j)^2\right\}\{1+o(1)\}}{3n(1+k-j)^2h(t_j,t_k)\left\{a(t_0)+b(t_0)\right\}^2(\log n)^2}\,,
\end{align*}
with an analogous upper bound for the terms, involving $Q_{k,j,}$, with $k<j$,
and, finally, if $k<j$, and $t_k,t_j\in[\e,1-\e]$,
\begin{align*}
&Ew_{j,k}^2\left\{\frac{f_0(t_k)\int_{I_k}\left(t-t_k\right)\,d\left(\G_{n,2}-G_2\right)(t)}{\G_{n,2}(t_{k+1})-\G_{n,2}(t_k)}\right\}^2
1_{\{M_k>0\}}\\
&\le\frac{3f_0(t_k)^2\{1+o(1)\}}{n(1+k-j)^2g_2(t_k)\left\{a(t_0)+b(t_0)\right\}^2(\log n)^2}\,.
\end{align*}
The terms for $t_k>1-\e$ are treated by using
\begin{align*}
&E w_{j,k}^2\left\{\frac{f_0(t_k)\int_{I_k}\left(t-t_k\right)\,d\G_{n,1}(t)}{\G_{n,1}(t_{k+1})-\G_{n,1}(t_k)}\right\}^2
\le E w_{j,k}^2f_0(t_k)^2\left(t_{k+1}-t_k\right)^2\\
&\le\frac{w_{j,k}^2f_0(t_k)^2}{K^2}
\sim \frac{9a(t_k)^2f_0(t_k)^2}{K^2(k-j+1)^2\left\{a(t_0)+b(t_0)\right\}^2(\log n)^2}\\
&\qquad\qquad\qquad\,\sim\frac{9a(t_k)^2f_0(t_k)^2}{K^4(t_k-t_j)^2\left\{a(t_0)+b(t_0)\right\}^2(\log n)^2}\,,
\end{align*}
with a similar upper bound for $t_k<\e$ and
$$
E w_{j,k}^2\left\{\frac{f_0(t_k)\int_{I_k}\left(t-t_k\right)\,d\G_{n,2}(t)}{\G_{n,2}(t_{k+1})-\G_{n,2}(t_k)}\right\}^2.
$$

We also have, for example, if $k'>k>j$,
\begin{align*}
&E\left\{\frac{w_{j,k}f_0(t_k)\int_{I_k}\left(t-t_k\right)\,d\left(\G_{n,1}-G_1\right)(t)}{\G_{n,1}(t_{k+1})-\G_{n,1}(t_k)}\right.\\
&\left.\qquad\qquad\qquad\cdot
\frac{w_{j,k'}f_0(t_{k'})\int_{I_{k'}}\left(t-t_{k'}\right)\,d\left(\G_{n,1}-G_1\right)(t)}{\G_{n,1}(t_{k'+1})-\G_{n,1}(t_{k'})}\right\}\\
&\sim\frac{9a(t_k)a(t_{k'})}{4nK^2\{a(t_0)+b(t_0)\}^2(k-j+1)(k'-j+1)(\log n)^2},
\end{align*}
and the expectation of other cross-product terms can be treated similarly.

Combining these results, we find that the variance of the conditional expectation
$$
\a_n^{-1}E\left\{\widetilde{F}_n(t_j)-F_0(t_j)\bigm|N_k,Q_{j,k},\,k>j;\,M_k,Q_{k,j},\,k<j\right\}
$$
is of order $O(1/\{n(\log n)^2\a_n^2\})=o(1)$.\\
(ii). We have, if $t_k\in[\e,1-\e]$,
\begin{align*}
\frac{\int_{I_k}\left\{F_0(t)-F_0(t_k)\right\}\,d\G_{n,1}(t)}{\G_{n,1}(t_{k+1})-\G_{n,1}(t_k)}
&=\frac{\tfrac12c^2f_0(t_k)\left(t_{k+1}-t_k\right)^2g_1(t_k)\left\{1+o_p(1)\right\}}
{cg_1(t_k)\left(t_{k+1}-t_k\right)\left\{1+o_p(1)\right\}}\\
&=\tfrac12cf_0(t_k)\left(t_{k+1}-t_k\right)\left\{1+o_p(1)\right\},
\end{align*}
and similarly,
$$
\frac{\int_{I_k}\left\{F_0(t)-F_0(t_k)\right\}\,d\G_{n,1}(t)}{\G_{n,2}(t_{k+1})-\G_{n,2}(t_k)}
=\tfrac12cf_0(t_k)\left(t_{k+1}-t_k\right)\left\{1+o_p(1)\right\},
$$
Moreover, if $k>j$,
\begin{align*}
&\frac{\int_{t\in I_j,\,u\in I_k}\left\{F_0(u)-F_0(t)-F_0(t_k)+F_0(t_j)\right\}\,d\H_n(t,u)}
{\int_{t\in I_j,\,u\in I_k}\,d\H_n(t,u)}\\
&=\tfrac12c\frac{\left\{f_0(t_k)\left(t_{k+1}-t_k\right)-f_0(t_j)\left(t_{j+1}-t_j\right)\right\}h(t_j,t_k)\left\{1+o_p(1)\right\}}
{h(t_j,t_k)\left\{1+o_p(1)\right\}}\\
&=\tfrac12c\left\{f_0(t_k)\left(t_{k+1}-t_k\right)-f_0(t_j)\left(t_{j+1}-t_j\right)\right\}\left\{1+o_p(1)\right\},
\end{align*}
with a similar expansion for $k<j$. The $o_p(1)$-terms are uniform in $k$, as follows by using exponential inequalities of the same type as used in Lemma \ref{lemma:asy}.

It is easily seen that the terms, involving values of $t_k\notin[\e,1-\e]$ give a negligible contribution, by noting that
$$
\frac{f_0(t_k)\int_{I_k}\left(t-t_k\right)\,d\G_{n,1}(t)}{\G_{n,1}(t_{k+1})-\G_{n,1}(t_k)}\le f_0(t_k)\left(t_{k+1}-t_k\right),
$$
if $k>j$, with a similar upper bound if $t_k<t_j$.
The results now follows by multiplying with $w_{j,k}$ and summing over $k$, see (\ref{heuristic_expansion}).
\eop

\vspace{0.3cm}
We now define
\begin{equation}
\label{U_{n,k}}
U_{n,k}=n^{-1}\left\{N_k'-E\left\{N_k'|N_k\right\}\right\},
\end{equation}
and
\begin{equation}
\label{V_{n,k}}
V_{n,k}=n^{-1}\left\{M_k'-E\left\{M_k'|M_k\right\}\right\}.
\end{equation}
Note that these are the numerators of the ``variance parts" in (\ref{condit_exp_decomp1}) and (\ref{condit_exp_decomp2}), divided by $n$.
The following lemma shows that (in the proper scaling for Birg\'e's statistic) the variances of the sums of terms, involving $U_{n,k}$
and $V_{n,k}$ in Birg\'e's statistic, tend to zero.

\begin{lemma}
\label{lemma:negligeability}
Let the conditions of Theorem \ref{limdis} be satisfied, let $t_j=t_0$, and let $\a_n$ be defined by (\ref{alpha_n}).
Moreover, let $U_{n,k}$ and $V_{n,k}$ be defined by (\ref{U_{n,k}}) and (\ref{V_{n,k}}).
Then, as $n\to\infty$,
\begin{align*}
&\a_n^{-2}\mbox{\rm var}\left(\sum_{k:k>j}\frac{w_{j,k}U_{n,k}}{\G_{n,1}\left(t_{k+1}\right)-\G_{n,1}\left(t_k\right)}
+\sum_{k:k< j}\frac{w_{j,k}V_{n,k}}{\G_{n,2}\left(t_{k+1}\right)-\G_{n,2}\left(t_k\right)}\right)\to0.
\end{align*}
\end{lemma}

\noindent
{\bf Proof.} We have:
\begin{align*}
&\mbox{\rm var}\left(\sum_{k:k>j}\frac{w_{j,k}U_{n,k}}{\G_{n,1}\left(t_{k+1}\right)-\G_{n,1}\left(t_k\right)}
+\sum_{k:k< j}\frac{w_{j,k}V_{n,k}}{\G_{n,1}\left(t_{k+1}\right)-\G_{n,1}\left(t_k\right)}\right)\\
&=\sum_{k:k>j}\mbox{\rm var}\left(\frac{w_{j,k}U_{n,k}}{\G_{n,1}\left(t_{k+1}\right)-\G_{n,1}\left(t_k\right)}\right)
+\sum_{k:k<j}\mbox{\rm var}\left(\frac{w_{j,k}V_{n,k}}{\G_{n,2}\left(t_{k+1}\right)-\G_{n,2}\left(t_k\right)}\right),
\end{align*}
since the covariances of the terms in the sum are zero. As before, we define the ratios to be zero if the denominator is zero.

Furthermore:
$$
\mbox{\rm var}\left(\frac{w_{j,k}U_{n,k}}{\G_{n,1}\left(t_{k+1}\right)-\G_{n,1}\left(t_k\right)}\right)
=E\left(\frac{w_{j,k}U_{n,k}}{\G_{n,1}\left(t_{k+1}\right)-\G_{n,1}\left(t_k\right)}\right)^2,
$$
since $E\,w_{j,k}U_{n,k}/\{\G_{n,1}\left(t_{k+1}\right)-\G_{n,1}\left(t_k\right)\}=0$.
Noting that the weights $w_{j,k}$ have upper bound
$$
\frac{\sqrt{n\left\{\G_{n,1}\left(t_{k+1}\right)-\G_{n,1}\left(t_k\right)\right\}}}{(k-j+1)W_j}\,,
$$
we now obtain:
\begin{align*}
&\a_n^{-2}\mbox{\rm var}\left(\frac{w_{j,k}U_{n,k}}{\G_{n,1}\left(t_{k+1}\right)-\G_{n,1}\left(t_k\right)}\right)\\
&\le\a_n^{-2}E\frac{nU_{n,k}^2}{(k-j+1)^2\left\{\G_{n,1}\left(t_{k+1}\right)-\G_{n,1}\left(t_k\right)\right\}W_j^2}\\
&=\a_n^{-2}E\frac{\int_{I_k}F_0(t)\left\{1-F_0(t)\right\}\,d\G_{n,1}(t)}{(k-j+1)^2\left\{\G_{n,1}\left(t_{k+1}\right)-\G_{n,1}\left(t_k\right)\right\}
W_j^2}\\
&=\frac{\a_n^{-2}\left\{F_0(t_k)\left(1-F_0(t_k)\right)+o(1)\right\}}{(k-j+1)^2}E\left\{1/W_j^2\right\}1_{\{W_j>0\}}\,.
\end{align*}
where (as before),
$$
U_{n,k}/\left\{\G_{n,1}\left(t_{k+1}\right)-\G_{n,1}\left(t_k\right)\right\}\stackrel{\mbox{def}}=0,
$$
if
$\G_{n,1}\left(t_{k+1}\right)-\G_{n,1}\left(t_k\right)=0$.

By (\ref{E1/W_j^m}):
\begin{equation}
\label{E1/W_j^2b}
E\left\{1/W_j^2\right\}1_{\{W_j>0\}}\sim\frac{K}{n\left\{a(t_0)+b(t_0)\right\}^2(\log n)^2}\asymp
n^{-2/3}(\log n)^{-5/3},
\end{equation}
So we obtain
\begin{align*}
&\a_n^{-2}\sum_{k:k>j}\mbox{\rm var}\left(\frac{
w_{j,k}U_{n,k}}{\G_{n,1}\left(t_{k+1}\right)-\G_{n,1}\left(t_k\right)}\right)\\
&\qquad\qquad+
\a_n^{-2}\sum_{k:k< j}\mbox{\rm
var}\left(\frac{ w_{j,k}V_{n,k}}{\G_{n,2}\left(t_{k+1}\right)-\G_{n,2}\left(t_k\right)}\right)\\
&=O\left(1/\log n\right).
\end{align*}
\eop

\vspace{0.3cm}
We now define, if $j<k$,
\begin{equation}
\label{W_{n,j,k}}
W_{n,j,k}=(\log n)\left\{Q_{j,k}'-E\left\{Q_{j,k}'|Q_{j,k}\right\}\right\},
\end{equation}
and, if $j>k$:
\begin{equation}
\label{W_{n,j,k}2}
W_{n,j,k}=(\log n)\left\{Q_{k,j}'-E\left\{Q_{k,j}'|Q_{k,j}\right\}\right\}\,.
\end{equation} 
Lemma \ref{lemma:negligeability} suggests that if $(n\log n)^{1/3}\{\widetilde{F}_n(t_0)-F_0(t_0)\}$ has a nondegenerate distribution,
this has to come from the sum:
\begin{equation}
\label{centering2}
-\sum_{k:k>j}w_{j,k}\frac{W_{n,j,k}}{c^2h(t_j,t_k)}
+\sum_{k:k<j}w_{j,k}\frac{W_{n,j,k}}{c^2h(t_k,t_j)}\,,
\end{equation}
The following lemma shows that (\ref{centering2}), with the random weights $w_{j,k}$ replaced by the deterministic weights
$\widetilde{w}_{j,k}$ indeed has a nondegenerate limit distribution.

\begin{lemma}
\label{lemma:lim_nondegen}
Let the conditions of Theorem \ref{limdis} be satisfied, let $t_j=t_0$. Moreover, let $W_{n,j,k}$
be defined by (\ref{W_{n,j,k}}) and (\ref{W_{n,j,k}2}). Then:
$$
-\sum_{k:k>j}\widetilde{w}_{j,k}\frac{W_{n,j,k}}{c^2h(t_j,t_k)}
+\sum_{k:k<j}\widetilde{w}_{j,k}\frac{W_{n,j,k}}{c^2h(t_k,t_j)}
\stackrel{{\cal D}}{\longrightarrow}
N\left(0,\s_0^2\right)
$$
where the right-hand side
denotes a normal random variable, with expectation $0$ and variance $\s_0^2$, defined by (\ref{sigma_0}) in Theorem \ref{limdis}.
\end{lemma}

\noindent
{\bf Proof.} We will prove the result by constructing a martingale-difference array, and applying Theorem 1,
p.\ 171 of \cite{Pol:84}. Define, for $k>j$, the random variables
$$
\xi_{n,k}=-\widetilde{w}_{j,k}\frac{W_{n,j,k}}{c^2h(t_j,t_{k})}.
$$
For $k<j$ we define
$$
\xi_{n,k}=\widetilde{w}_{j,k}\frac{W_{n,j,k}}{c^2h(t_k,t_j)},
$$
and (for notational convenience) we define $\xi_{n,j}\equiv0$.

Let the increasing sequence of $\s$-fields ${\cal F}_{n,k}$, $k=0,1,\dots$ be defined by
$$
{\cal F}_{n,0}=\emptyset,\,{\cal F}_{n,k}=\s\left\{\left(T_i,U_i,\dd_i\right),\,T_i\le t_{k+1},\,U_i\in I_j\right\},\,k\le j,
$$
and
$$
{\cal F}_{n,k}=\s\left\{\left(T_i,U_i,\dd_i\right),\,T_i\in I_j,\,U_i\le t_{k+1}\right\},\,j<k,
$$
where $\dd_i=(\dd_{i,1},\dd_{i,2},\dd_{i,3})$, as before. Note: $I_k=[t_k,t_{k+1})$, $k<K$, and $I_K=[t_K,t_{K+1}]$, under scheme (i),
and $I_k=[t_k,t_{k+1})$, $k\le K$, and $I_{K+1}=[t_{K+1},t_{K+2}]$ under scheme (ii) at the beginning of this section.

Then:
\begin{equation}
\label{martcond1}
E\left\{\xi_{n,k}\bigm|{\cal F}_{n,k-1}\right\}=0,\,k=1,2,\dots
\end{equation}
Here and in the following the indices $k$ run from 1 to $K$ or to $K+1$, depending on whether scheme (i) or (ii) holds, respectively.

Note that, if $k<j$, we can write
\begin{align*}
W_{n,j,k}&=n\log n\int_{(t,u)\in I_k\times I_j}\left\{\d_2-\left\{F_0(u)-F_0(t)\right\}\right\}\,d\P_n(t,u,\d)\\
&=\log n\sum_{i=1}^n\left\{\dd_{2,i}-\left\{F_0(u)-F_0(t)\right\}\right\}1_{\{T_i\in I_k,\,U_i\in I_j\}}\,.
\end{align*}
and that
$$
E\left\{\dd_{2,i}-\left\{F_0(U_i)-F_0(T_i)\right\}\bigm|{\cal F}_{n,k}\right\}=0,
$$
if $t_{k}<T_i<t_j$ and $U_i\in I_j$, using the independence of the $X_i$ from the pairs $(T_i,U_i)$. Similar relations hold
if $t_i\in I_j$. This implies
\begin{equation}
\label{martcond2}
E\left\{\xi_{n,k}\bigm|{\cal F}_{n,k-1}\right\}=0,\,k=1,2,\dots\,.
\end{equation}
It is also clear that $\xi_{n,k}$ is measurable with respect to ${\cal F}_{n,k}$.

Let the conditional variances $v_{n,k}$ be defined by
$$
v_{n,k}=E\left\{\xi_{n,k}^2\bigm|{\cal F}_{n,k-1}\right\},\,k=1,2,\dots\,.
$$
We first consider the indices $k$ such that
$$
|j-k|<\e_n K,
$$
where $\epsilon_n=(\log n)^{-1/3}$.
We then get, if $k<j$,
\begin{align*}
&v_{n,k}\\
&=\frac{\widetilde{w}_{j,k}^2n(\log n)^2}{c^4h(t_k,t_0)^2}\\
&\qquad\cdot E\left\{\int_{t\in I_k,\,u\in
I_j}\left\{F_0(u)-F_0(t)\right\}\left\{1-F_0(u)+F_0(t)\right\}\,d\H_n(t,u)
\bigm|{\cal F}_{n,k-1}\right\}\\
&=\frac{\widetilde{w}_{j,k}^2n^{1/3}(\log n)^{4/3}(t_0-t_k)f_0(t_0)\left\{1+o_p(1)\right\}}{c^2h(t_k,t_0)}\\
&=\frac{9\b(t_0)^2(n\log n)^{1/3}(j-k+1)f_0(t_0)\left\{1+o_p(1)\right\}}{c^2K h(t_k,t_0)(j-k+1)^2\log n}\\
&=\frac{9b(t_0)^2f_0(t_0)\left\{1+o_p(1)\right\}}{c\left\{a(t_0)+b(t_0)\right\}^2h(t_k,t_0)(j-k+1)\log n}\,.
\end{align*}
We similarly get:
\begin{align*}
v_{n,k}=\frac{9a(t_0)^2f_0(t_0)\left\{1+o_p(1)\right\}}{c\left\{a(t_0)+b(t_0)\right\}^2h(t_k,t_0)(j-k+1)\log n}\,.
\end{align*}
if $k>j$ and $k-j<\e_n K$.
The terms $v_{n,k}$, where $|k-j|\ge\e_n K$, give a negligible contribution, since
\begin{align*}
\sum_{k:j-k\ge\e_n K}v_{n,k}&=
\sum_{k:j-k\ge\e_n K}\frac{\widetilde{w}_{j,k}^2n(\log n)^2}{c^4h(t_k,t_0)^2} O_p\left\{(n\log n)^{-2/3}\right\}\\
&=O_p\left((\log n)^{-2/3}\right),
\end{align*}
using
$$
\sum_{k:j-k\ge\e_n K}\widetilde{w}_{j,k}^2=O\left(\sum_{k:j-k\ge\e_n K}\frac1{(j-k)^2(\log
n)^2}\right)=O\left(n^{-1/3}(\log n)^{-2}\right),
$$
as $n\to\infty$. So we find
\begin{equation}
\label{asym_condit_var}
\sum_{k:k\ne j}v_{n,k}\stackrel{p}\longrightarrow\s_0^2\,,
\end{equation}
since
$$
\sum_{m:m<\e_n K}\frac1{m+1}\sim\frac13\log n,\,n\to\infty.
$$

To get asymptotic normality, it only remains to show that the Lindeberg-type condition
\begin{equation}
\label{Lindeberg_condition}
\sum_{k\ne j} E\left\{\xi_{n,k}^2 1_{\left\{\left|\xi_{n,k}\right|>\e\right\}}\bigm|{\cal
F}_{n,k-1}\right\}\stackrel{p}\longrightarrow0,
\end{equation}
holds for each $\e>0$, since in that case both conditions of Theorem 1 of \cite{Pol:84} are satisfied.
To this end we use the conditional Cauchy-Schwarz inequality
\begin{equation}
\label{Cauchy-Schwarz}
E\left\{\xi_{n,k}^2 1_{\left\{\left|\xi_{n,k}\right|>\e\right\}}\bigm|{\cal
F}_{n,k-1}\right\}\le
\sqrt{E\left\{\xi_{n,k}^4\bigm|{\cal F}_{n,k-1}\right\}E\left\{1_{\left\{\left|\xi_{n,k}\right|>\e\right\}}\bigm|{\cal
F}_{n,k-1}\right\}}\,.
\end{equation}
Note that:
\begin{align}
\label{CS-second_term}&E\left\{1_{\left\{\left|\xi_{n,k}\right|>\e\right\}}\bigm|{\cal F}_{n,k-1}\right\}
\le\e^{-2}E\left\{\xi_{n,k}^2\bigm|{\cal F}_{n,k-1}\right\}
=\e^{-2}v_{n,k}=O_p(1/\log n),
\end{align}
as $n\to\infty$.
Using again the conditional independence of the $X_i$, given the values of the pairs $(T_i,U_i)$, and defining
$p_0(t,u)=F_0(u)-F_0(t)$, $\overline{p}_0(t,u)=1-p_0(t,u)$, we get, if $k<j$:
\begin{align}
\label{tightness_sum}
&E\left\{\xi_{n,k}^4\bigm|{\cal F}_{n,k-1}\right\}\nonumber\\
&\sim\frac{\widetilde{w}_{j,k}^4n(\log n)^4}{c^8h(t_k,t_0)^4}\cdot\nonumber\\
&\quad\cdot E\left\{\int_{t\in I_k,\,u\in
I_j}p_0(t,u)\overline{p}_0(t,u)\}\left\{p_0(t,u)^3+\overline{p}_0(t,u)^3\right\}\,d\H_n(t,u)
\bigm|{\cal F}_{n,k-1}\right\}\nonumber\\
&\qquad+\frac{\widetilde{w}_{j,k}^4n^2(\log n)^4}{c^8h(t_k,t_0)^4}\cdot\nonumber\\
&\qquad\quad\cdot E\left\{\biggl\{\int_{t\in I_k,\,u\in I_j}p_0(t,u)\{1-\overline{p}_0(t,u)\}\,d\H_n(t,u)\biggr\}^2\biggm|{\cal F}_{n,k-1}\right\}.
\end{align}
The first conditional expectation on the right-hand side of (\ref{tightness_sum}) arises from terms of the form
$$
E\left\{\left\{\dd_{2,i}-\left(F_0(U_i)-F_0(T_i)\right)\right\}^4\bigm|{\cal F}_{n,k-1}\right\},
$$
where $T_i\in I_k,\,U_i\in I_j$, and the second one from terms of the form
$$
E\left\{\left\{\dd_{2,i}-\left(F_0(U_i)-F_0(T_i)\right)\right\}^2
\left\{\dd_{2,i'}-\left(F_0(U_{i'})-F_0(T_{i'})\right)\right\}^2\bigm|{\cal
F}_{n,k-1}\right\},
$$
where $i\ne i'$ and $T_i,T_{i'}\in I_k;\,U_i,U_{i'}\in I_j$, where we added the diagonal terms (where $i=i'$) for simplicity of
notation, since they give a negligible contribution. The other conditional expectations of crossproducts are zero. If $k>j$ we get an
entirely similar expansion, with the roles of $t$ and
$u$ interchanged.

The first term on the right-hand side of (\ref{tightness_sum}) gives a contribution of order $O_p(1/\sqrt{\log n})$ in the summation of
the terms
$$
\sqrt{E\left\{\xi_{n,k}^4\bigm|{\cal F}_{n,k-1}\right\}}
$$
over $k$. The square root of the second term is of order $O_p(1/\{|j-k|\log n\})$, if $|j-k|<\e_n K$, which leads to a contribution of
order $O_p(1)$ in the above summation. The part where $|j-k|\ge\e_n K$ is again negligible.

So we get, using (\ref{Cauchy-Schwarz}) and (\ref{CS-second_term}),
\begin{align*}
&\sum_{k=1}^K E\left\{\xi_{n,k}^2 1_{\left\{\left|\xi_{n,k}\right|>\e\right\}}\bigm|{\cal F}_{n,k-1}\right\}
=O_p\left(1/\sqrt{\log n}\right)\sum_{k=1}^K\sqrt{E\left\{\xi_{n,k}^4\bigm|{\cal F}_{n,k-1}\right\}}\\
&=O_p\left(1/\sqrt{\log n}\right).
\end{align*}
\eop

\vspace{0.4cm}
\noindent
{\bf Proof of Theorem \ref{limdis}.}\\
ad (i). Lemma \ref{lemma:negligeability} shows that the terms involving $N_k'/N_k$ and $M_k'/M_k$ only give a contribution to
the asymptotic bias of Birg\'e's statistic, but not to the limit distribution of the centered part. The limit distribution of the centered
part therefore arises from the terms $W_{n,j,k}$, where
$$
W_{n,j,k}=n\log n\int_{(t,u)\in I_k\times I_j}\left\{\d_2-\left\{F_0(u)-F_0(t)\right\}\right\}\,d\P_n(t,u,\d),
$$
which are the numerators of the fractions
\begin{align*}
&\frac{(n\log n)^{1/3}\left\{Q_{j,k}'-E\left(Q_{j,k}'|Q_{j,k}\right)\right\}}{Q_{j,k}}\\
&=\frac{n\log n\int_{(t,u)\in I_k\times I_j}\left\{\d_2-\left\{F_0(u)-F_0(t)\right\}\right\}\,d\P_n(t,u,\d)}
{(n\log n)^{2/3}\int_{(t,u)\in I_j\times I_k}d\H_n(t,u)}.
\end{align*}
Now note that
$$
(n\log n)^{2/3}\int_{(t,u)\in I_j\times I_k}d\H_n(t,u)=c^2h(t_j,t_k)\left\{1+o_p(1)\right\}.
$$
where the $o_p(1)$-term is uniform in $k$ by the results, given above. Moreover, by part (i) of Lemma
\ref{lemma:asy},
\begin{align*}
\sum_{k\ne j}w_{j,k}\frac{W_{n,j,k}}{c^2\widetilde{h}(t_j,t_k)}=\sum_{k\ne j}\widetilde{w}_{j,k}\frac{W_{n,j,k}}{c^2\widetilde{h}(t_j,t_k)}
+o_p(1)\sum_{k\ne j}\frac{W_{n,j,k}}{(k-j+1)\log n}\,,
\end{align*}
where
$$
\widetilde{h}(t,u)=h(t,u),\,t<u,\,\widetilde{h}(t,u)=h(u,t),\,t\ge u.
$$
The result now follows from Lemma \ref{lemma:lim_nondegen}.\\
ad (ii).
We first prove (\ref{bias_lim}). Since $E\widetilde{F}_n(t_j)$ is the expectation of
$$
E\left\{\widetilde{F}_n(t_j)\bigm|N_k,Q_{j,k},\,k>j;\,M_k,Q_{k,j},\,k<j\right\},
$$
part (i) of Lemma \ref{lemma:bias} tells us that
$$
\a_n^{-2}E\left\{E\left\{\widetilde{F}_n(t_j)\bigm|N_k,Q_{j,k},\,k>j;\,M_k,Q_{k,j},\,k<j\right\}-E\widetilde{F}_n(t_j)\right\}^2\to0,
$$
as $n\to\infty$.
This implies:
$$
\a_n^{-1}\left\{E\left\{\widetilde{F}_n(t_j)\bigm|N_k,Q_{j,k},\,k>j;\,M_k,Q_{k,j},\,k<j\right\}-E\widetilde{F}_n(t_j)\right\}
\stackrel{p}\longrightarrow0,
$$
as $n\to\infty$. But since, by part (ii) of Lemma \ref{lemma:bias},
$$
\a_n^{-1}\left\{E\left\{\widetilde{F}_n(t_j)\bigm|N_k,Q_{j,k},\,k>j;\,M_k,Q_{k,j},\,k<j\right\}-F_0(t_0)\right\}
\stackrel{p}\longrightarrow\tfrac12cf_0(t_0),
$$
as $n\to\infty$, we must have:
$$
\a_n^{-1}\left\{E\widetilde{F}_n(t_j)-F_0(t_0)\right\}\to\tfrac12cf_0(t_0),\,n\to\infty.
$$
This yields (\ref{bias_lim}).

To prove (\ref{var_lim}), we first note that, by part (i) of Lemma \ref{lemma:bias}, the variance of the conditional expectation
$$
\a_n^{-1}E\left\{\widetilde{F}_n(t_j)-F_0(t_0)\bigm|N_k,Q_{j,k},\,k>j;\,M_k,Q_{k,j},\,k<j\right\}
$$
in the decomposition
\begin{align*}
&\a_n^{-1}\left\{\widetilde{F}_n(t_j)-F_0(t_0)\right\}\\
&=\a_n^{-1}\left\{\widetilde{F}_n(t_j)-E\left\{\widetilde{F}_n(t_j)\bigm|N_k,Q_{j,k},\,k>j;\,M_k,Q_{k,j},\,k<j\right\}\right\}\\
&\qquad\qquad+\a_n^{-1}E\left\{\widetilde{F}_n(t_j)-F_0(t_0)\bigm|N_k,Q_{j,k},\,k>j;\,M_k,Q_{k,j},\,k<j\right\}
\end{align*}
tends to zero. By Lemma \ref{lemma:negligeability} the sum of terms involving $N_k$ and $M_k$ also gives an asymptotically negligible 
contribution to $\a_n^{-1}\left\{\widetilde{F}_n(t_j)-F_0(t_0)\right\}$.

So we only have to consider the contribution of terms of the form
\begin{equation}
\label{Q_{j,k}_conditional1}
\frac{\a_n^{-1}w_{j,k}\left\{Q_{j,k}'-E\left(Q_{j,k}'|Q_{j,k}\right)\right\}}{Q_{j,k}},\,k>j,
\end{equation}
and
\begin{equation}
\label{Q_{j,k}_conditional2}
\frac{\a_n^{-1}w_{j,k}\left\{Q_{k,j}'-E\left(Q_{k,j}'|Q_{k,j}\right)\right\}}{Q_{k,j}},\,k<j.
\end{equation}
The variance of (\ref{Q_{j,k}_conditional1}) is given by
$$
E\frac{n(\log n)^2w_{j,k}^2\int_{(t,u)\in I_k\times I_j}\left\{F_0(u)-F_0(t)\right\}\left\{1-\left(F_0(u)-F_0(t)\right)\right\}\,d\H_n(t,u)}
{(n\log n)^{4/3}\left\{\int_{(t,u)\in I_j\times I_k}d\H_n(t,u)\right\}^2}\,.
$$
Lemma \ref{lemma:asy} gives (uniform) exponential inequalities are derived for the probabilities of the events of the following type:
$$
A_{j,k}\stackrel{\mbox{\small def}}=\left\{\left|(n\log n)^{2/3}\int_{(t,u)\in I_j\times I_k}d\H_n(t,u)-c^2h(t_j,t_k)\right|>\e
c^2h(t_j,t_k)\right\},
$$
yielding upper bounds, tending to zero faster than any power of $n$. So we get:
\begin{align*}
&E\Biggl\{\frac{n(\log n)^2w_{j,k}^2\int_{(t,u)\in I_k\times I_j}\left\{F_0(u)-F_0(t)\right\}\left\{1-\left(F_0(u)-F_0(t)\right)\right\}\,d\H_n(t,u)}
{(n\log n)^{4/3}\left\{\int_{(t,u)\in I_j\times I_k}d\H_n(t,u)\right\}^2}\\
&\qquad\qquad\qquad\qquad\qquad\qquad\qquad\qquad\qquad\qquad\qquad\qquad\qquad\qquad\qquad\cdot 1_{A_{j,k}}\Biggr\}\\
&\le K\left\{F_0(t_{k+1})-F_0(t_j)\right\}E\frac{(n\log n)^{2/3}\{1+o(1)\}} {(1+k-j)^2W_j^2}1_{A_{j,k}\cap\{W_j>0\}}\\
&\le K^3\left\{F_0(t_{k+1})-F_0(t_j)\right\}E\frac{(n\log n)^{2/3}\{1+o(1)\}} {(1+k-j)^2} P\left(A_{j,k}\right),
\end{align*}
which tends to zero faster than any power of $n$, uniformly in $k$. Here we use the lower bound $1/K$ for $W_j1_{\{W_j>0\}}$.

So we find:
\begin{align*}
&E\frac{n(\log n)^2w_{j,k}^2\int_{(t,u)\in I_k\times
I_j}\left\{F_0(u)-F_0(t)\right\}(\left\{1-\left(F_0(u)-F_0(t)\right)\right\}\,d\H_n(t,u)} {(1-\e)^2c^4h(t_j,t_k)}\\
&\qquad\qquad\qquad\qquad\qquad\qquad\qquad\qquad\qquad\qquad\qquad\qquad\qquad\qquad+o(1/K)\\
&\ge E\frac{n(\log n)^2w_{j,k}^2\int_{(t,u)\in I_k\times
I_j}\left\{F_0(u)-F_0(t)\right\}\left\{1-\left(F_0(u)-F_0(t)\right)\right\}\,d\H_n(t,u)} {(n\log n)^{4/3}\left\{\int_{(t,u)\in I_j\times
I_k}d\H_n(t,u)\right\}^2}\\
&\ge E\frac{n(\log n)^2w_{j,k}^2\int_{(t,u)\in I_k\times
I_j}\left\{F_0(u)-F_0(t)\right\}\left\{1-\left(F_0(u)-F_0(t)\right)\right\}\,d\H_n(t,u)} {(1+\e)^2c^4h(t_j,t_k)}\\
&\qquad\qquad\qquad\qquad\qquad\qquad\qquad\qquad\qquad\qquad\qquad\qquad\qquad\qquad+o(1/K).
\end{align*}
This implies:
\begin{align*}
&E\frac{n(\log n)^2w_{j,k}^2\int_{(t,u)\in I_k\times
I_j}\left\{F_0(u)-F_0(t)\right\}\left\{1-\left(F_0(u)-F_0(t)\right)\right\}\,d\H_n(t,u)} {(n\log n)^{4/3}\left\{\int_{(t,u)\in I_j\times
I_k}d\H_n(t,u)\right\}^2}\\
&= E\frac{n(\log n)^2 w_{j,k}^2\int_{(t,u)\in I_k\times
I_j}\left\{F_0(u)-F_0(t)\right\}\left\{1-\left(F_0(u)-F_0(t)\right)\right\}\,d\H_n(t,u)} {c^4h(t_j,t_k)}\\
&\qquad\qquad\qquad\qquad\qquad\qquad\qquad\qquad\qquad\qquad\qquad\qquad\qquad\qquad+o(1/K).
\end{align*}

Now let, for $t_k<1-\d$, where $\d>0$, the event $B_k$ be defined by
$$
B_k\stackrel{\mbox{\small def}}=\left\{(1+k-j)\left|(n\log n)^{1/3}\int_{u\in I_k}\,d\G_n(u)-cg_1(t_k)\right|>\e
cg_1(t_k)\right\},
$$
For $t_k\ge1-\d$, we define the event $B_k$ by:
$$
B_k\stackrel{\mbox{\small def}}=\left\{(1+k-j)\left|(n\log n)^{1/3}\int_{u\in I_k}\,d\G_n(u)-cg_1(t_k)\right|>\e
c\right\},
$$
Similarly to what is true for $A_{j,k}$, we have that $P(B_k)$ tends to zero faster than any power of $n$, uniformly in $k$. 
This shows that 
we also can replace $w_{j,k}$ by $\widetilde{w}_{j,k}$ in the asymptotic expression for the variance, using the fact that the terms for
$t_k>1-\d$ will give a contribution of lower order in the summation. So we find:
\begin{align*}
&n(\log n)^2\cdot\\
&\cdot\sum_{k:j>k}\frac{ E w_{j,k}^2\int_{(t,u)\in I_k\times
I_j}\left\{F_0(u)-F_0(t)\right\}\left\{1-\left(F_0(u)-F_0(t)\right)\right\}d\H_n(t,u)} {c^4h(t_j,t_k)}\\
&\sim\sum_{k:j>k}\frac{9na(t_k)^2E\int_{(t,u)\in I_k\times
I_j}\left\{F_0(u)-F_0(t)\right\}\left\{1-\left(F_0(u)-F_0(t)\right)\right\}d\H_n(t,u)}
{c^4\left\{a(t_0)+b(t_0)\right\}^2(j-k+1)^2h(t_j,t_k)}\\
&\sim\sum_{k:j>k}\frac{9na(t_k)^2
\left\{F_0(t_k)-F_0(t_j)\right\}\left\{1-\left(F_0(t_k)-F_0(t_j)\right)\right\}(n\log n)^{-2/3}}
{c^2\left\{a(t_0)+b(t_0)\right\}^2(j-k+1)^2h(t_j,t_k)}\\
&\sim\sum_{k:j>k}\frac{9na(t_j)^2
f_0(t_j)\left(t_k-t_j\right)(n\log n)^{-2/3}}
{c^2\left\{a(t_0)+b(t_0)\right\}^2(j-k+1)^2h(t_j,t_k)}\\
&\sim\sum_{k:j>k}\frac{9na(t_j)^2 f_0(t_j)(n\log n)^{-1}}{c\left\{a(t_0)+b(t_0)\right\}^2(j-k+1)h(t_j,t_k)}
\sim\frac{3a(t_0)^2 f_0(t_0)}{c\left\{a(t_0)+b(t_0)\right\}^2h(t_j,t_k)}.
\end{align*}
Similarly we find that the summation for $k<j$ gives a contribution which is asymptotically equivalent to
$$
\frac{3b(t_0)^2 f_0(t_0)}{c\left\{a(t_0)+b(t_0)\right\}^2h(t_j,t_k)}.
$$
This yields (\ref{var_lim}).\eop

\noindent
{\bf Proof of Lemma \ref{lemma:asy}}.\\
We first prove (\ref{W_j-asymp}).
By Bennett's inequality (see, e.g., \cite{Pol:84}, p.\ 192) we have, for $\e>0$,
\begin{align*}
&P\left\{\left|N_k/n-EN_k/n\right|>\frac{\e}{K}\right\}\\
&\le 2\exp\left\{-\frac{n\e^2}{2K^2\int_{t\in I_k} g_1(t)\,dt}\,
\f\left(\frac{\e}{K\int_{t\in I_k}g_1(t)\,dt}\right)\right\},
\end{align*}
where
\begin{equation}
\label{K-L-function}
\f(x)=\frac{2\left\{(1+x)\log(1+x)-x\right\}}{x^2}\,,\,x>0.
\end{equation}
This way of stating Bennett's inequality first appeared in \cite{galen:80}.
The function $\f$ satisfies $\lim_{x\downarrow0}\f(x)=1$ and
$$
\f(x)\ge \frac1{1+x/3}\,,\,x>0,
$$
see \cite{Pol:84}, p.\ 192, p.\ 193.

By the continuity of $g_1$ on $[0,1]$ there exists for each $k$ a $\xi_k\in I_k$ such that
$$
\int_{I_k}g_1(t)\,dt=g_1(\xi_k)\left\{t_{k+1}-t_k\right\}.
$$
Hence we get, for each $k$,
\begin{align*}
&P\left\{\left|N_k/n-EN_k/n\right|>\frac{\e}{K}\right\}
\le 2\exp\left\{-\frac{n\e^2}{2Kg_1(\xi_k)}\,
\f\left(\frac{\e}{g_1(\xi_k)}\right)\right\}\\
&=2\exp\left\{-\frac{cn^{2/3}\e^2}{2g_1(\xi_k)(\log n)^{1/3}}\,
\f\left(\frac{\e}{g_1(\xi_k)}\right)\right\}.
\end{align*}
Similarly we get, for each $k$ and points $\eta_k\in I_k$,
\begin{align*}
&P\left\{\left|M_k/n-EM_k/n\right|>\frac{\e}{K}\right\}
\le2\exp\left\{-\frac{cn^{2/3}\e^2}{2g_2(\eta_k)(\log n)^{1/3}}\,
\f\left(\frac{\e}{2g_2(\eta_k)}\right)\right\}.
\end{align*}
Moreover, if $j<k$,
\begin{align*}
&P\left\{\left|Q_{j,k}/n-EQ_{j,k}/n\right|>\frac{\e}{K^2}\right\}\\
&\le 2\exp\left\{-\frac{n\e^2}{2K^4\int_{t\in I_j,\,u\in I_k} h(t,u)\,dt\,du}\,
\f\left(\frac{\e}{K^2\int_{t\in I_j,\,u\in I_k} h(t,u)\,dt\,du}\right)\right\}\\
&=2\exp\left\{-\frac{c^2n^{1/3}\e^2\{1+o(1)\}}{2h(t_j,t_k)(\log n)^{2/3}}\,
\f\left(\frac{\e\{1+o(1)\}}{2h(t_j,t_k)}\right)\right\}.
\end{align*}
with a similar upper bound, if $k<j$.

Let $\e>0$, let $\overline{h}$ be defined by
\begin{equation}
\label{overline{h}}
\overline{h}(t,u)=h(t,u),\,u\ge t,\,\overline{h}(t,u)=h(u,t),\,u<t,
\end{equation}
and similarly $\overline{Q}_{j,k}$ by
\begin{equation}
\label{overline{Q}}
\overline{Q}_{j,k}(t,u)=Q_{j,k}(t,u),\,u\ge t,\,k\ge j,\,\overline{Q}_{j,k}(t,u)=Q_{k,j}(u,t),\,(u,t),\,u<t,\,k<j.
\end{equation}
Moreover, let the set $A_{j,\e}$ be defined by
\begin{align*}
A_{j,\e}&=\left\{\sup_{k\ne
j}\left|\overline{Q}_{j,k}/n-E\overline{Q}_{j,k}/n\right|\le\frac{\e}{K^2},\,\sup_{k>j}\left|N_k/n-EN_k/n\right|\le\frac{\e}{K},\right.\\
&\left.\qquad\qquad\qquad\qquad\qquad\qquad\qquad\qquad\qquad\qquad\sup_{k<j}\left|M_k/n-EM_k\right|\le\frac{\e}{K}\right\}.
\end{align*}
and let
\begin{equation}
\label{h_j}
h_j=\inf_{k:k\ne j}\overline{h}(t_j,t_k).
\end{equation}

Then we have:
\begin{equation}
\label{W_j-upperbound}
1-P\left(A_{j,\e}\right)=O\left(n^{1/3}\exp\left\{-\frac{c^2n^{1/3}\e^2}{4h_j(\log n)^{2/3}}\,
\f\left(\frac{\e}{4h_j}\right)\right\}\right).
\end{equation}
Furthermore, as $n\to\infty$,
$$
\sup_{k:k>j}\left|KEN_k/n-g_1(t_k)\right|=\sup_{k:k>j}\left|K\int_{t\in I_k}g_1(t)\,dt-g_1(t_k)\right|\to0,
$$
also on the last interval, since $g_1(t)\to0$ on this interval,
$$
\sup_{k:k<j}\left|KEM_k/n-g_2(t_k)\right|=\sup_{k:k<j}\left|K\int_{t\in I_k}g_2(t)\,dt-g_2(t_k)\right|\to0,
$$
also on the first interval, since $g_2(t)\to0$ on this interval.
We also have:
$$
\left|K^2E\overline{Q}_{j,k}/n-\overline{h}(t_j,t_k)\right|=\left|K^2\int_{t\in I_j,\,u\in I_k}
\overline{h}(t,u)\,dt\,du-\overline{h}(t_j,t_k)\right|\to0,
$$
uniformly for all $t_k$, not belonging to the first or last interval, which may not have length $1/K$ (see the construction of the
intervals of Birg\'e's statistic at the beginning of section \ref{sec:Birge}). But on these intervals we have
$$
h(t,t_j)\wedge g_2(t)=g_2(t)\mbox{ and }h(t_j,t)\wedge g_1(t)=g_1(t),
$$
respectively. So we get:
\begin{equation}
\label{right-upperbound}
\sup_{k:k>j}\left|(KEN_k/n)\wedge \left(K^2EQ_{j,k}/n\right)-g_1(t_k)\wedge h(t_j,t_k)\right|\to0,
\end{equation}
and
\begin{equation}
\label{left-upperbound}
\sup_{k:k<j}\left|(KEM_k/n)\wedge \left(K^2EQ_{k,j}/n\right)-g_2(t_k)\wedge h(t_k,t_j)\right|\to0,
\end{equation}
Hence, we get from (\ref{W_j-upperbound}), (\ref{right-upperbound}) and (\ref{left-upperbound}), on the set $A_{j,\e}$,
\begin{align*}
W_j&=\sum_{k<j}\frac{\sqrt{M_k\wedge(KQ_{k,j})}}{j-k+1}+\sum_{k>j}\frac{\sqrt{N_k\wedge(KQ_{j,k})}}{k-j+1}\\
&=\sqrt{n}\sum_{k<j}\frac{\sqrt{(M_k/n\wedge(KQ_{k,j}/n)}}{j-k+1}+\sqrt{n}\sum_{k>j}\frac{\sqrt{(N_k/n)\wedge(KQ_{j,k}/n)}}{k-j+1}\\
&\ge\sqrt{n(1-\e)}\left\{\sum_{k<j}\frac{\sqrt{(EM_k/n)\wedge(KEQ_{k,j}/n)}}{j-k+1}\right.\\
&\left.\qquad\qquad\qquad\qquad\qquad\qquad\qquad\qquad+\sum_{k>j}\frac{\sqrt{(EN_k/n)\wedge(KEQ_{j,k}/n)}}{k-j+1}\right\}
\end{align*}
\begin{align*}
&=\frac{\sqrt{n(1-\e)}}{\sqrt{K}}\left\{\sum_{k<j}\frac{\sqrt{(KEM_k/n)\wedge(K^2EQ_{k,j}/n)}}{j-k+1}\right.\\
&\left.\qquad\qquad\qquad\qquad\qquad\qquad+\sum_{k>j}\frac{\sqrt{(KEN_k/n)\wedge(K^2EQ_{j,k}/n)}}{k-j+1}\right\}\\
&=\frac{\sqrt{n(1-\e)}}{\sqrt{K}}\left\{\sum_{k<j}\frac{\sqrt{g_2(t_k)\wedge h(t_k,t_j)}}{j-k+1}
+\sum_{k>j}\frac{\sqrt{g_1(t_k)\wedge h(t_j,t_k)}}{k-j+1}\right\}\,,
\end{align*}
and similarly
$$
W_j\le\frac{\sqrt{n(1+\e)}}{\sqrt{K}}\left\{\sum_{k<j}\frac{\sqrt{g_2(t_k)\wedge h(t_k,t_j)}}{j-k+1}
+\sum_{k>j}\frac{\sqrt{g_1(t_k)\wedge h(t_j,t_k)}}{k-j+1}\right\}\,.
$$
Moreover, letting $\e_n=(\log n)^{-1/3}$, we get:,
\begin{align*}
&\sum_{k<j}\frac{\sqrt{g_2(t_k)\wedge h(t_k,t_j)}}{j-k+1}
+\sum_{k>j}\frac{\sqrt{g_1(t_k)\wedge h(t_j,t_k)}}{k-j+1}\\
&=\sum_{k:t_j-t_k<\e_n}\frac{\sqrt{g_2(t_k)\wedge h(t_k,t_j)}}{j-k+1}
+\sum_{k:t_k-t_j<\e_n}\frac{\sqrt{g_1(t_k)\wedge h(t_j,t_k)}}{k-j+1}\\
&\qquad+\sum_{k:t_j-t_k\ge\e_n}\frac{\sqrt{g_2(t_k)\wedge h(t_k,t_j)}}{j-k+1}
+\sum_{k:t_k-t_j\ge\e_n}\frac{\sqrt{g_1(t_k)\wedge h(t_j,t_k)}}{k-j+1}\\
&=\tfrac13\left\{a(t_0)+b(t_0)\right\}(\log n)\{1+o(1)\}.
\end{align*}
Relation (\ref{W_j-asymp}) now follows.

To prove (\ref{E1/W_j^m}) we first note that
$$
E\frac1{W_j^m}1_{\{W_j>0\}\cap A_{j,\e}^c}=O\left((K+1)^m n^{1/3}\exp\left\{-\frac{c^2n^{1/3}\e^2}{4h_j(\log n)^{2/3}}\,
\f\left(\frac{\e}{4h_j}\right)\right\}\right),
$$
where $h_j$ is defined by (\ref{h_j}), since $W_j\ge1/(K+1)$, if $W_j>0$
Thus we find:
\begin{align*}
&E\frac1{W_j^m}1_{\{W_j>0\}}=E\frac1{W_j^m}1_{\{W_j>0\}\cap A_{j,\e}}+E\frac1{W_j^m}1_{\{W_j>0\}\cap A_{j,\e}^c}\\
&\le\frac1{(1-\e)^{m/2}}\left\{\sum_{k<j}\frac{\sqrt{EM_k\wedge(K EQ_{k,j})}}{j-k+1}+\sum_{k>j}\frac{\sqrt{EN_k\wedge(K
EQ_{j,k})}}{j-k+1}\right\}^{-m}\\
&\qquad+O\left((K+1)^m n^{1/3}\exp\left\{-\frac{c^2n^{1/3}\e^2}{4h_j(\log n)^{2/3}}\,
\f\left(\frac{\e}{4h_j}\right)\right\}\right)\\
&\sim \left(\frac{9K}{n(1-\e)}\right)^{m/2}\left\{\left(a(t_0)+b(t_0)\right)\log n\right\}^{-m}\,,\,n\to\infty.
\end{align*}
and similarly
$$
E\frac1{W_j^m}1_{\{W_j>0\}}\ge\left(\frac{9K}{n(1+\e)}\right)^{m/2}\left\{\left(a(t_0)+b(t_0)\right)\log n\right\}^{-m}\,,\,n\to\infty,
$$
implying
$$
E\frac1{W_j^m}1_{\{W_j>0\}}\sim\left(\frac{9K}{n}\right)^{m/2}\left\{\left(a(t_0)+b(t_0)\right)\log n\right\}^{-m}\,,\,n\to\infty,
$$
which proves (\ref{E1/W_j^m}).

Finally we get for $j>k$:
\begin{align*}
&(1+k-j)E\left|w_{j,k}-\widetilde{w}_{j,k}\right|1_{\{W_j>0\}}\\
&=E\left|\frac{\sqrt{N_k\wedge(KQ_{k,j})}}{W_j}1_{\{W_j>0\}}-\frac{3a(t_k)}{\left\{a(t_0)+b(t_0)\right\}\log n}\right|\\
&\le E\left|\frac{\sqrt{(KN_k/n)\wedge(K^2Q_{k,j}/n)}-3a(t_k)}{W_j\sqrt{K/n}}\right|1_{\{W_j>0\}}\\
&\qquad\qquad+3a(t_k)E\left|\frac{\sqrt{n/K}}{W_j}1_{\{W_j>0\}}-\frac1{\left\{a(t_0)+b(t_0)\right\}\log n}\right|.
\end{align*}
Applying the Cauchy-Schwarz inequality on the first term on the right-hand side yields, if $j<k$,
\begin{align*}
&E\left|\frac{\sqrt{(KN_k/n)\wedge(K^2Q_{j,k}/n)}-3a(t_k)}{W_j\sqrt{K/n}}\right|1_{\{W_j>0\}}\\
&\le\left\{E\left\{\sqrt{(KN_k/n)\wedge(K^2Q_{j,k}/n)}-3a(t_k)\right\}^2\right\}^{1/2}\sqrt{E1/W_j^2}\\
&=o(1/\log n),
\end{align*}
uniformly in $k$, using (\ref{E1/W_j^m}) and the exponential inequalities for
$$
P\left\{\left|N_k/n-EN_k/n\right|>\frac{\e}{K}\right\}\mbox{ and }P\left\{\left|Q_{j,k}/n-EQ_{j,k}/n\right|>\frac{\e}{K^2}\right\}
$$
derived above. Using ((\ref{E1/W_j^m})) again, we get that the second term satisfies the inequality
\begin{align*}
&3a(t_k)E\left|\frac{\sqrt{n/K}}{W_j}1_{\{W_j>0\}}-\frac1{\left\{a(t_0)+b(t_0)\right\}\log n}\right|\\
&\le\left\{E\left\{\frac{\sqrt{n/K}}{W_j}1_{\{W_j>0\}}-\frac1{\left\{a(t_0)+b(t_0)\right\}\log n}\right\}^2\right\}^{1/2}
=o(1/\log n).
\end{align*}
The case $k<j$ is treated similarly.

We also have:
\begin{align*}
&(1+k-j)E\left|w_{j,k}-\widetilde{w}_{j,k}\right|1_{\{W_j=0\}}
=(1+k-j)\left|\widetilde{w}_{j,k}\right|P\left\{W_j=0\right\}\\
&=o(1/\log n),
\end{align*}
since, in fact, $P\left\{W_j=0\right\}$ tends to zero exponentially fast in $n$. This proves (\ref{uniformity}).
\eop

\vspace{0.4cm}
We next discuss the proof of Theorem \ref{limdis_sep}.
Since the following lemmas have proofs analogous to the proofs of Lemma \ref{lemma:bias} and Lemma \ref{lemma:negligeability}
in section \ref{sec:Birge} we omit their proofs.

\begin{lemma}
\label{lemma:deterministic_denom_sep}
Let the observation density $h$, the number of intervals $K$ and the constant $c$ be
as in Theorem \ref{limdis}, and let $t_k=t_k^{(n)}$ be the left boundary point of a sub-interval of the partition, Moreover, let $F_0$
have a continuous derivative on $[0,1]$, and let $\G_{n,1}$ and $\G_{n,2}$ be the empirical distribution functions of the $T_i$ and
$U_i$, respectively. Then
\begin{enumerate}
\item[(i)]
\begin{equation}
\label{first_term_sep}
\left\{\frac{N_k'}{N_k}-F_0(t_k)\right\}1_{\{N_k>0\}}=
\frac{U_{n,k}+\int_{t\in I_k}\left\{F_0(t)-F_0(t_k)\right\}\,d\G_{n,1}(t)}
{\G_{n,1}\left(t_{k+1}\right)-\G_{n,1}\left(t_k\right)}\,,
\end{equation}
where
$$
U_{n,k}=\int_{t\in I_k}\left\{\d_1-F_0(t)\right\}\,d\P_n(t,u,\d).
$$
\item[(ii)]
\begin{equation}
\label{second_term_sep}
\left\{\frac{M_k'}{M_k}-F_0(t_k)\right\}1_{\{M_k>0\}}=
\frac{V_{n,k}+\int_{t\in I_k}\left\{F_0(t)-F_0(t_k)\right\}\,d\G_{n,2}(t)}
{\G_{n,2}\left(t_{k+1}\right)-\G_{n,2}\left(t_k\right)}\,,
\end{equation}
where
$$
V_{n,k}=\int_{u\in I_k}\left\{\d_1+\d_2-F_0(u)\right\}\,d\P_n(t,u,\d).
$$
\item[(iii)]
Let $t_j=t_0$. Then, if $k>j$,
\begin{align}
\label{third_term_sep}
&n^{1/3}\left\{\frac{Q_{j,k}'}{Q_{j,k}}-\left\{F_0(t_k)-F_0(t_j)\right\}\right\}\nonumber\\
&=\left\{\tfrac12c\left\{f_0(t_k)-f_0(t_0)\right\}+\frac{W_{n,j,k}}{c^2h(t_j,t_k)}\right\}\left\{1+o_p(1)\right\}\,,
\end{align}
where
$$
W_{n,j,k}=n\int_{(t,u)\in I_j\times I_k}\left\{\d_2-\left\{F_0(u)-F_0(t)\right\}\right\}\,d\P_n(t,u,\d).
$$
If $k<j$ we get:
\begin{align}
&n^{1/3}\left\{\frac{Q_{j,k}'}{Q_{j,k}}-\left\{F_0(t_0)-F_0(t_k)\right\}\right\}\nonumber\\
&=\left\{\frac{W_{n,j,k}}{c^2h(t_k,t_j)}+\tfrac12c\left\{f_0(t_0)-f_0(t_k)\right\}\right\}\left\{1+o_p(1)\right\},
\end{align}
where
$$
W_{n,j,k}=n\int_{(t,u)\in I_k\times I_j}\left\{\d_2-\left\{F_0(u)-F_0(t)\right\}\right\}\,d\P_n(t,u,\d).
$$
\item[(iv)] The $o_p(1)$ terms in (iii) tend to zero uniformly in $k$.
\end{enumerate}
\end{lemma}

\begin{lemma}
\label{lemma:negligeability_sep}
Let the conditions of Theorem \ref{limdis_sep} be satisfied, and let $t_j=t_0$.
Then, as $n\to\infty$,
\begin{align*}
&n^{2/3}\mbox{\rm var}\left(\sum_{k:k>j}\frac{w_{j,k}U_{n,k}}{\G_{n,1}\left(t_{k+1}\right)-\G_{n,1}\left(t_k\right)}
+\sum_{k:k< j}\frac{w_{j,k}V_{n,k}}{\G_{n,2}\left(t_{k+1}\right)-\G_{n,2}\left(t_k\right)}\right)\to0.
\end{align*}
\end{lemma}

Since the first moment of the asymptotic distribution follows in a similar way as in section \ref{sec:Birge}, using the
representations of the components $N_k'/N_k$, etc.\ of Lemma \ref{lemma:deterministic_denom_sep}, the proof of Theorem \ref{limdis_sep}
again boils down to proving the following lemma.

\begin{lemma}
\label{lemma:lim_nondegen_sep}
Let the conditions of Theorem \ref{limdis_sep} be satisfied, and let $t_j=t_0$. Moreover, let $W_{n,j,k}$
be defined as in part (iii) of Lemma \ref{lemma:deterministic_denom_sep}, and let $\s^2$ be defined as in Theorem \ref{limdis_sep}. Then:
$$
-\sum_{k:k>j}\widetilde{w}_{j,k}\frac{W_{n,j,k}}{c^2h(t_j,t_k)}
+\sum_{k:k<j}\widetilde{w}_{j,k}\frac{W_{n,j,k}}{c^2h(t_k,t_j)}
\stackrel{{\cal D}}{\longrightarrow}
N\left(0,\s^2\right),
$$
where the right-hand side
denotes a normal random variable, with expectation $0$ and variance $\s^2$, defined by (\ref{as_var_sep}).
\end{lemma}

\noindent
{\bf Proof.} Since the proof follows the same lines as the proof of Theorem \ref{limdis}, we only give the main steps. We define
the martingale difference array in the same way as in the proof of Theorem \ref{limdis}. Then, if $k<j$, we get the following
representation of  the conditional variance
\begin{align*}
&v_{n,k}\\
&=\frac{n\widetilde{w}_{j,k}^2}{c^4h(t_k,t_0)^2}\\
&\qquad\cdot E\left\{\int_{t\in I_k,\,u\in
I_j}\left\{F_0(u)-F_0(t)\right\}\left\{1-F_0(u)+F_0(t)\right\}\,d\H_n(t,u)
\bigm|{\cal F}_{n,k-1}\right\}\\
&=\frac{n^{1/3}\widetilde{w}_{j,k}^2\left\{F_0(t_0)-F_0(t_k)\right\}\left\{1-F_0(t_0)+F_0(t_k)\right\}\left\{1+o_p(1)\right\}}{c^2h(t_k,t_0)}.
\end{align*}
Similarly we get, if $k<j$,
\begin{align*}
&v_{n,k}\\
&=\frac{n\widetilde{w}_{j,k}^2}{c^4h(t_k,t_0)^2}\\
&\qquad\cdot E\left\{\int_{t\in I_j,\,u\in
I_k}\left\{F_0(u)-F_0(t)\right\}\left\{1-F_0(u)+F_0(t)\right\}\,d\H_n(t,u)
\bigm|{\cal F}_{n,k-1}\right\}\\
&=\frac{n^{1/3}\widetilde{w}_{j,k}^2\left\{F_0(t_k)-F_0(t_0)\right\}\left\{1-F_0(t_k)+F_0(t_0)\right\}\left\{1+o_p(1)\right\}}
{c^2h(t_0,t_k)}.
\end{align*}
Hence, using (\ref{modified_w_sep}) and a Riemann sum approximation, we obtain:
\begin{equation}
\label{condit_var_conv}
\sum_{k\ne j}v_{n,k}\stackrel{p}\longrightarrow \s^2.
\end{equation}

It remains to show  the Lindeberg-type condition
\begin{equation}
\label{Lindeberg_condition_sep}
\sum_{k\ne j} E\left\{\xi_{n,k}^2 1_{\left\{\left|\xi_{n,k}\right|>\d\right\}}\bigm|{\cal
F}_{n,k-1}\right\}\stackrel{p}\longrightarrow0,
\end{equation}
for each $\d>0$.
We again use the conditional Cauchy-Schwarz inequality
\begin{equation}
\label{Cauchy-Schwarz_sep}
E\left\{\xi_{n,k}^2 1_{\left\{\left|\xi_{n,k}\right|>\d\right\}}\bigm|{\cal
F}_{n,k-1}\right\}\le
\sqrt{E\left\{\xi_{n,k}^4\bigm|{\cal F}_{n,k-1}\right\}E\left\{1_{\left\{\left|\xi_{n,k}\right|>\d\right\}}\bigm|{\cal
F}_{n,k-1}\right\}}\,.
\end{equation}
We have:
\begin{align}
\label{CS-second_term_sep}
E\left\{1_{\left\{\left|\xi_{n,k}\right|>\d\right\}}\bigm|{\cal F}_{n,k-1}\right\}
\le\frac1{\d^2}E\left\{\xi_{n,k}^2\bigm|{\cal F}_{n,k-1}\right\}
=O_p\left(1/K\right)=O_p\left(n^{-1/3}\right).
\end{align}
Letting
$p_0(t,u)=F_0(u)-F_0(t)$, $\overline{p}_0(t,u)=1-p_0(t,u)$, we get, if $k<j$:
\begin{align}
\label{tightness_sum_sep}
&E\left\{\xi_{n,k}^4\bigm|{\cal F}_{n,k-1}\right\}\nonumber\\
&\sim\frac{n\widetilde{w}_{j,k}^4}{c^8h(t_k,t_0)^4}\nonumber\\
&\qquad\cdot E\left\{\int_{t\in I_k,\,u\in
I_j}p_0(t,u)\overline{p}_0(t,u)\}\left\{p_0(t,u)^3+\overline{p}_0(t,u)^3\right\}\,d\H_n(t,u)
\bigm|{\cal F}_{n,k-1}\right\}\nonumber\\
&\qquad+\frac{n^2\widetilde{w}_{j,k}^4}{c^8h(t_k,t_0)^4}
E\left\{\biggl\{\int_{t\in I_k,\,u\in I_j}p_0(t,u)\{1-\overline{p}_0(t,u)\}\,d\H_n(t,u)\biggr\}^2\biggm|{\cal F}_{n,k-1}\right\}.
\end{align}
The first and second terms on the right-hand side are, respectively, of order
$$
O_p\left(\frac1{K^3(t_0-t_k)^4}\right)\mbox{ and }O_p\left(\frac1{K^2(t_0-t_k)^4}\right).
$$
So the second term is dominant, and hence:
\begin{align}
\label{tightness_sum_sep2}
&\sum_{k<j}\frac{n\widetilde{w}_{j,k}^2}{c^4h(t_k,t_0)^2}\nonumber\\
&\cdot\left\{E\left\{\biggl\{\int_{t\in I_k,\,u\in I_j}p_0(t,u)
\{1-\overline{p}_0(t,u)\}\,d\H_n(t,u)\biggr\}^2\biggm|{\cal F}_{n,k-1}\right\}\right\}^{1/2}\nonumber\\
&=O_p\left(\frac1{K}\sum_{k<j}\frac1{(t_0-t_k)^2}\right)=O_p\left(\int_\e^{t_0-\e}\frac1{(t_0-t)^2}\,dt\right)
=O_p(1).
\end{align}
Similarly the sum of the terms for $k>j$ is $O_p(1)$. The result now follows from (\ref{Cauchy-Schwarz_sep}) and
(\ref{CS-second_term_sep}).\eop

\vspace{0.3cm}
The proof of Theorem \ref{limdis_sep} can now be finished by making the transition from the random weights to the deterministic
weights, using  Lemma \ref{lemma:asy2} (see the proof of Theorem \ref{limdis} at the end of section \ref{sec:Birge}), and using the
central limit result of Lemma \ref{lemma:lim_nondegen_sep}.

\bibliographystyle{amsplain}
\bibliography{references}

\end{document}